\documentclass[titlepage,11pt]{article}
\usepackage{latexsym}
\usepackage{amsfonts}
\usepackage{amsmath}
\newcommand\blackslug{\hbox{\hskip 1pt \vrule width 4pt height 8pt depth 1.5pt
        \hskip 1pt}}
\newcommand\bbox{\hfill \quad \blackslug \bigbreak}
\newcommand{\vare}{\varepsilon}
\newcommand{\ins}{\operatorname{ins}}
\newcommand{\bd}{\operatorname{bd}}
\newcommand{\mac}{\mathcal}
\def\DD{\hbox{-}}
\def\CC{\hbox{-}\cdots\hbox{-}}
\def\LL{,\ldots,}
\def\cupcup{\cup\cdots\cup}

\title{Excluding disjoint Kuratowski graphs}
\author{Neil Robertson\\
Ohio State University, Columbus, OH 43210
\\
\\
Paul Seymour\thanks{Supported by AFOSR grant
FA9550-22-1-0234, and NSF grant  DMS-2154169.}\\
Princeton University, Princeton, NJ 08544}

\date{January 10, 2024; revised \today}

\newtheorem{thm}{}[section]

\newcommand{\Proof}{\noindent{\bf Proof.}\ \ }

\begin{document}
\maketitle
\begin{abstract}
A graph is a {\em $k$-Kuratowski graph} if it has exactly $k$ components, each isomorphic to $K_5$ or to $K_{3,3}$.
We prove that if a graph $G$ contains no $k$-Kuratowski graph as a minor,
then there is a set $X$ of boundedly many vertices such that $G\setminus X$
can be drawn in a (possibly disconnected) surface in which no $k$-Kuratowski graph can be drawn.

\end{abstract}

\section{Introduction}
Graphs in this paper are finite, and have no loops or parallel edges. The Kuratowski-Wagner~\cite{kuratowski,
wagner} theorem says that:
\begin{thm}\label{KWthm}
A graph contains neither of $K_5, K_{3,3}$ as a minor (or equivalently, has no subgraph which is a subdivision
of $K_5$ or $K_{3,3}$) if and only if it is planar.
\end{thm}
There is a similar theorem that characterizes the minimal graphs that cannot be drawn on the projective 
plane~\cite{archdeacon}, but for higher surfaces it seems hopeless to obtain the complete list of 
excluded minors (although we know that the list is finite for every surface, by~\cite{GM8,GM20}). Even for the torus,
the list numbers in the tens of thousands (at least 17,523 according to Myrvold and Woodcock~\cite{myrvold}).

We need a few definitions. 
Let us say a graph is a {\em $k$-Kuratowski graph} if it has exactly $k$ components, each isomorphic to $K_5$ or to $K_{3,3}$.
A {\it surface} $\Sigma$ is a non-null compact $2$-manifold, with (possibly null)
boundary, and possibly disconnected, and $\hat{\Sigma}$ is the surface without boundary obtained from $\Sigma$ by pasting a closed disc onto each component of its boundary.
(The presence of the boundary makes no different to which graphs can be drawn in the surface, but we will need the boundary later.)
If $\Sigma$ is a connected surface, and $\Sigma$ is orientable, its {\em genus} is the number of handles we add to a 2-sphere to 
make $\hat{\Sigma}$,
and if $\Sigma$ is non-orientable, 
its {\em genus} is the number of crosscaps we add to a 2-sphere to make $\hat{\Sigma}$.  The genus of 
a general (disconnected) surface is the sum of the genera of its components. This definition is non-standard, but 
convenient for our purposes; because there is a theorem that if $H$ is  a $k$-Kuratowski graph and $\Sigma$ is a surface, then $H$ can be drawn in $\Sigma$ 
if and only if the genus of $\Sigma$ is at least $k$.

Thus, excluding the $(k+1)$-Kuratowski graphs as minors is necessary for embeddability
in a surface of genus at most $k$. Unsurprisingly, the converse is false, 
but our main theorem says that it is not {\em very} false:
\begin{thm}\label{thm1}
For every integer $k\ge 0$, there is a number $f(k)$ with the following property. If $G$ is a graph with no $(k+1)$-Kuratowski graph
as a minor, 
then there exists $X\subseteq V(G)$ with $|X|\le f(k)$ 
such that $G\setminus X$ can be drawn in a surface of genus at most $k$.
\end{thm}
The proof of \ref{thm1} is an application of the graph minors structure theorem of~\cite{GM16}. 
We first proved 
it in the early 1990's, but did not write it up at that time.
There was another, related  result we proved at that time and did not write up, that will appear in a subsequent paper~\cite{Ksums}.
It says, roughly, that if a graph does not contain a $k$-Kuratowski graph as a minor, and does not contain as a minor the graph
made from $k$ copies of $K_5$ or $K_{3,3}$ by identifying together one vertex of each, and also does not contain similar graphs made
from $k$ copies of $K_5$ and $K_{3,3}$ by identifying pairs or triples of vertices, then $G$ has bounded genus.

Before we go on,
let us prove the result mentioned earlier, that:
\begin{thm}\label{embedding}
For all integers $k\ge 0$, if $\Sigma$ is a surface and $H$ is a $k$-Kuratowski graph, then $H$ can be drawn in $\Sigma$ if and only if $\Sigma$ has genus at least $k$.
\end{thm}
\Proof
(Thanks to Carsten Thomassen and Bojan Mohar for their help with this proof.)
Certainly $H$ can be drawn in any surface with genus at least $k$, and we prove the converse by induction on the sum of $k$ and 
the number of components of $\Sigma$. Let $\ell$ be the genus of $\Sigma$.
We are assuming that $H$ can be drawn in $\Sigma$, and need to show that $\ell\ge k$.
Since $\Sigma$ has genus at least zero, we may assume that $k\ge 1$.
We may assume that $\Sigma$ has null boundary. Suppose that $\Sigma$ is disconnected, and is the disjoint union of two surfaces
$\Sigma_1,\Sigma_2$. For $i = 1,2$, let $H_i$ be the subgraph of $H$ drawn in $\Sigma_i$; then $H_i$ is a $k_i$-Kuratowski graph,
for some $k_i$, where $k_1+k_2=k$. From the inductive hypothesis $\Sigma_i$ has genus at least $k_i$ for $i = 1,2$, and so the 
result holds. Hence we may assume that $\Sigma$ is connected. 

The orientable genus of a graph equals the sum of the orientable genera of its components, by a theorem of 
Battle et al.~\cite{battle}; so $H$ has orientable genus at least $k$. 
The {\em Euler genus} of a connected surface $\Sigma$ is twice the genus of $\Sigma$ if $\Sigma$ is orientable, and equal 
to its genus otherwise; and the {\em Euler genus} of a graph is the minimum Euler genus of the connected surfaces in which it can be drawn.
A theorem of Miller~\cite{miller} implies that the Euler genus of any graph is the sum of the Euler genera 
of its components, 
and it follows that $H$ has Euler genus $k$. So if $\Sigma$ is orientable, then $\ell \ge k$ since $H$ has orientable genus at least $k$; and if $\Sigma$ is non-orientable, then $\ell\ge k$ since $H$ has Euler genus at least $k$.
This proves \ref{embedding}.~\bbox

\section{A lemma about packing and covering}

We will need a result from~\cite{ding}, the following:

\begin{thm}\label{ding}
If $H$ is a hypergraph with every hyperedge nonempty, then 
$$\tau(H)\le 11\lambda(H)^2(\lambda(H)+\nu(H)+3)\binom{\lambda(H)+\nu(H)}{\nu(H)}^2.$$
\end{thm}
Here, a {\em hypergraph} means a finite set of subsets, called {\em hyperedges}, of a finite set of {\em vertices};
$\tau(H)$ means the size of the smallest set of vertices that has nonempty intersection with every hyperedge;
$\nu(H)$ means the maximum $k$ such that there are $k$ hyperedges, pairwise disjoint; and $\lambda(H)$ means the maximum $k$
such that there are hyperedges $A_1\LL A_k$ and vertices $v_{ij}\;(1\le i<j\le k)$ such that for $1\le i<j\le k$, $v_{ij}$
belongs to $A_i$ and to $A_j$ and no other other hyperedges among $A_1\LL A_k$.

Let us see how this will be applied. First, here is a result also proved in~\cite{ding}:
\begin{thm}\label{hitregions}
For every surface $\Sigma$ without boundary there is a function $f$ with the following property.
Let $G$ be a graph drawn on $\Sigma$; let $U\subseteq V(G)$ and let $S$ be a subset of the set of regions of $G$.
Then for every integer $k\ge 0$, either there are $k$ vertices in $U$ such that no region in $S$ is incident with any two of them,
or there is a set of at most $f(k)$ regions in $S$ such that every vertex in $U$ is incident with one of them.
\end{thm}
\Proof Let $H$ be the hypergraph with vertex set $S$, and with hyperedges the sets of all regions in $S$ incident with $v$,
for each $v\in U$. (Two of these hyperedges might be equal, and if so, omit one.) By \ref{ding}, it suffices to show that
$\lambda(H)$
is bounded.  Let $k=\lambda(H)$ and let $A_1\LL A_k$ be hyperedges and $v_{ij}\;(1\le i<j\le k)$ as in the definition of $\lambda(H)$.
For $1\le i\le k$, draw a vertex $v_i$ in each $r_i$, and take a line segment in $r_i$ between $v_i$ and each $v_{ij}\;(j>i)$
and each $v_{hi}\;(h<i)$,
all pairwise disjoint except for $v_i$. This makes a graph which is a subdivision of the complete graph on $k$ vertices.
Since the size of the largest complete graph that can be drawn in $\Sigma$ is bounded, it follows that $\lambda(H)$
is bounded. This proves \ref{hitregions}.~\bbox

Bienstock and Dean~\cite{bienstock} claimed to prove a much stronger result: that under the same hypotheses, the hypergraph $H$, 
constructed the same way, satisfies that $\nu(H)\ge c\tau(H)$, where the constant $c>0$ depends only on $\Sigma$. (Their proof
for higher surfaces has a mistake, which can easily be fixed, using 
that all graphs drawable on a fixed surface have bounded chromatic number, but we omit details since we will not use their result.)

We need another application of \ref{ding}. If $G$ is drawn in a surface $\Sigma$, let us say vertex $u$ 
{\em sees} vertex $v$ if there is a region incident with both $u,v$ (so $v$ sees itself). A {\em wheel neighbourhood}
of a vertex $v$ of $G$ means the set of all vertices that $v$ sees. We need:
\begin{thm}\label{hitwheels}
For every surface $\Sigma$, there is a function $f$ with the following property.
Let $G$ be a graph drawn on $\Sigma$, and let $U\subseteq V(G)$. 
Then for every integer $k\ge 0$, either there are $k$ vertices in $U$ such that their wheel neighbourhoods are pairwise disjoint,
or there is a set $X$ of at most $f(k)$ vertices of $G$ such that every $u\in U$ sees some vertex in $X$.
\end{thm}
\Proof Let $H$ be the hypergraph with vertex set $V(G)$, and with hyperedges the sets of all wheel neighbourhoods of vertices in $U$.
By \ref{ding}, it suffices to show that
$\lambda(H)$ is at most a constant depending on $\Sigma$. 
Let $k=\lambda(H)$;  and let $A_1\LL A_k$ be hyperedges and $v_{ij}=v_{ji}\;(1\le i<j\le k)$ as in the definition of 
$\lambda(H)$. For $1\le i\le k$, choose $v_i\in U$ such that $A_i$ is the wheel neighbourhood of $v_i$. 
For all distinct $i,j\in \{1\LL k\}$, let $r_{ij}$ be a region incident with $v_i$ and with $v_{ij}$. It follows that for all
$h\ne i,j$, $v_h$ is not incident with $r_{ij}$, since $v_h$ does not see $v_{ij}$. Hence there is a line segment between $v_i,v_j$, passing through $v_{ij}$ and with interior included in $r_{ij}\cup r_{ji}\cup \{v_{ij}\}$, and 
these line segments are pairwise disjoint except for their ends. Hence they make a drawing of a $k$-vertex complete graph in $\Sigma$, and so $\lambda(H)$ is bounded. This proves \ref{hitwheels}.~\bbox

\section{Paintings}

For brevity, let us say a {\em K-graph} in $G$ is a subgraph of $G$ that is isomorphic to a subdivision
of $K_5$ or $K_{3,3}$. We define the {\em K-number} of a graph $G$ to be the maximum number of K-graphs
in $G$ that are pairwise vertex-disjoint; that is, the maximum $k$ such that $G$ contains a $k$-Kuratowski graph as a minor.

The goal of this section and the next is to prove a very special case of \ref{thm1} (we will show later that the general problem
can be boiled down to this special case). The special case is essentially the following. Take a hypergraph, with hyperedges of 
size two or three, drawn on a surface without boundary; every hyperedge is represented by a closed disc, and its vertices are represented by nodes 
in the surface, drawn 
in the boundary of the disc. The discs for distinct hyperedges meet only in the nodes that represent vertices they have in common. Now make a graph 
by replacing each hyperedge of size three by a copy of $K_{2,3}$, where the three nodes of the hyperedge become the three vertices in one part of the bipartition of $K_{2,3}$. Replace each hyperedge of size two either by an edge, or a copy of $K_5\setminus e$
(that is, the graph obtained from $K_5$ by deleting one edge), or a copy of $K_{3,3}\setminus e$, where the two vertices of 
the hyperedge become the two ends of the deleted edge, making a graph $G$. We need to prove that $G$
satisfies \ref{thm1}. The difficulty is, the gadgets $K_{2,3}$, $K_5\setminus e$ and $K_{3,3}\setminus e$ that replace
the hyperedges are not K-graphs, and yet in general, $G$ is far from being embeddable in the surface.

So far we have beem somewhat casual about drawings in surfaces, but now we need to be more precise.
We need the definition of a painting from~\cite{GM17}, except for the moment 
we only need paintings in surfaces without boundary, which makes them considerably simpler. 
If $\Sigma$ is a surface and 
$X \subseteq \Sigma$, we denote the topological closure of $X$ by
$\overline{X}$, and we denote $\overline{X}\setminus X$ by
$\tilde{X}$.  A {\it painting} $\Gamma$ in a surface
$\Sigma$ without boundary is a pair $(U, N)$, where $U \subseteq \Sigma$
is closed and $N \subseteq U$ is finite, such that
\begin{itemize}
\item $U\setminus N$ has only finitely many arcwise
connected components (which we call {\it cells} of $\Gamma$); and
\item for each cell $c$, its closure
$\overline{c}$ is a disc and $\tilde{c}$ is a subset of the boundary
of this disc, and $|\tilde{c} | \le 3$. We say $c$ has {\em order $|\tilde{c}|$}. 
\end{itemize}
We define $U (\Gamma) = U$ and $N (\Gamma) = N$.  The members of
$N(\Gamma)$ are called the {\it nodes} of $\Gamma$.
A {\it region} of $\Gamma$
is an arcwise connected component of $\Sigma\setminus U$.  Thus, each region is a
connected open set. We define {\em incidence} between nodes, cells and regions as follows. 
If $n$ is a node, $c$ is a cell, and $r$ is a region, we say $c, n$ are incident if $n\in \tilde{c}$,
$r,n$ are incident if $n\in \overline{r}$, and $c,r$ are incident if $c\cap \overline{r}\ne \emptyset$. (We remind the reader that
$\tilde{c}$ is disjoint from $c$; a region with closure only intersecting $\tilde{c}$ is not incident with $c$.)
A subset $X \subseteq \Sigma$ is $\Gamma$-{\it normal}
if $X \cap U ( \Gamma ) \subseteq N ( \Gamma )$.
An {\em O-arc} in $\Sigma$ is a subset homeomorphic to a circle.
If $n,n'$ are nodes, we say $n$ {\em sees} $n'$ if there is a region incident with both $n,n'$. 

In this section we are given a painting $\Gamma$ in a surface $\Sigma$ without boundary, 
satisfying the following conditions:
\begin{itemize}
\item[{\bf (I1)}] $|\tilde{c} | \in \{2,3\}$ for every cell $c$.
\item[{\bf (I2)}] The closure of every region is a closed disc.
\item[{\bf (I3)}] If $r,r'$ are distinct regions, and $u,v\in N(\Gamma)$ are distinct and both incident with both $r,r'$,
then there is a cell $c$ with $\tilde{c}=\{u,v\}$, incident with both $r,r'$.
\item[{\bf (I4)}] For each cell $c$ of order two, let $r,r'$ be the regions incident with $c$; then
every node incident with both $r,r'$ is in $\tilde{c}$.
\item [{\bf (I5)}] For each cell $c$ of order three, and each node $n\notin \tilde{c}$, there exists $n'\in \tilde{c}$ such that 
no region is incident with both $n,n'$.
\item[{\bf (I6)}] For every $\Gamma$-normal O-arc $F$ with $|F\cap N(\Gamma)|\le 6$, there is a closed disc $\ins(F)\subseteq \Sigma$ with
boundary $F$ called the {\em inside} of $F$. (If $\Sigma$ is not a sphere, there is at most one choice for $\ins(F)$, but for a sphere, either of the discs bounded by $F$ might be its inside. But not both: one has been selected.) Let us say 
a {\em plate} of $\Gamma$ is a disc $\ins(F)$ for some $\Gamma$-normal O-arc $F$ with $|F\cap N(\Gamma)|\le 6$.

\end{itemize}

If $c$ is a cell of order two, then for each region $r$ incident with $c$, the closure of $r$ is called a {\em wing} of $c$.
If $c$ is a cell of order three, 
then for each node $n\in \tilde{c}$, the union of the closures of all regions incident with $n$ is called a {\em wing} of $c$. 
If $\mac B$ is a set of cells of $\Gamma$, every wing of a cell in $\mac B$ is called a {\em $\mac B$-wing}, and we call
the members of $\mac B$ {\em bodies}.
We will prove:
\begin{thm}\label{3planar}
For every surface $\Sigma$ without boundary,  there is a function $f$ such that, if $\Gamma$ satisfies {\bf (I1)}\LL {\bf (I6)},
 then for 
all $k$, and for all choices of the set $\mac B$ of bodies, 
either there are $k$ pairwise disjoint $\mac B$-wings, or there is a set $X$ of nodes with $|X|\le f(k)$, and a set $Y$ of plates with $|Y|\le f(k)$, 
such that for every body $c$, either $X\cap \tilde{c}\ne \emptyset$, or $c$ is a 
subset of a plate in $Y$.
\end{thm}

We will prove this in steps. First, we assume that every body is a cell of order two.

\begin{thm}\label{3planar2cell}
For every surface $\Sigma$ without boundary,  there is a function $f$ such that, if $\Gamma$ satisfies {\bf (I1)}\LL {\bf (I6)},
 then for
all $k$, and for all choices of the set $\mac B$ of bodies such that every body has order two,
either there are $k$ pairwise disjoint $\mac B$-wings, or there is a set $X$ of nodes with $|X|\le f(k)$, and a set $Y$ of plates with $|Y|\le f(k)$,
such that for every body $c$, either $X\cap \tilde{c}\ne \emptyset$, or $c$ is a
subset of a plate in $Y$.
\end{thm}
\Proof Let $f_1$ be the function $f$ of \ref{hitregions}, and let $f(k) = 3f_1(k)^2$. We will show that $f$ satisfies the theorem.
Let $k\ge 0$ be an integer. 
We assume that there do not exist $k$ pairwise disjoint $\mac B$-wings. From \ref{hitregions}, we deduce
that there exists $X\subseteq N(\Gamma)$ with $|X|\le f_1(k)$, such that $X\cap W\ne \emptyset$
for every wing $W$.

Let $x,x'\in X$ be distinct. We say a body $c$ is {\em pinned by} $\{x,x'\}$ if $x,x'\notin \tilde{c}$ and 
$x\in D$ and $x'\in D'$, where $D,D'$ are the wings of $c$.
\\
\\
(1) {\em For all distinct $x,x'\in X$, there is a set of at most three plates such that every body pinned by $\{x,x'\}$ is a subset of one of them.}
\\
\\
Let $\mathcal{C}$ be the set of all cells pinned by $\{x,x'\}$, and let 
$M$ be the union of $\tilde{c}$ over all $c\in \mac C$.
For each node $n\in M$, there is a region incident with both $x,n$ and not with $x'$; 
choose a line segment $L(n,x)$ between $n,x$ with interior
within one such region. Choose these 
all pairwise disjoint except for $x$.  Define $L(n,x')$ similarly. For each $n\in M$, $L(n,x)\cup L(n,x')$ is a $\Gamma$-normal
line segment between $x,x'$ that passes through $n$ and no other nodes except $x,x'$, and we denote it by $S_n$.
We call $S_n$ a {\em strut}.
We may assume that $|\mac C|\ge 2$, since otherwise $\mac C=\emptyset$ and the claim is true. 
If $n,n'\in M$ are distinct, then $S_{n}\cup S_{n'}$ is a $G$-normal O-arc passing though only four nodes, and 
so its inside is defined; let us denote it by $D_{n,n'}$. Choose $n,n'\in M$ such that $D(n,n')$ is maximal. If $M\subseteq D(n,n')$,
then every body pinned by $\{x,x'\}$ is a subset of the plate $D(n,n')$, except possibly for one cell incident with $n,n'$; 
let us take a second plate for this errant cell, and then the claim holds.

Thus we may assume that $M\not\subseteq D(n,n')$; choose $m\in M\setminus \{n,n'\}$. 
For every two of the three plates
$D(n,n'), D(m,n), D(m,n')$, either one includes the other, or their union is a closed disc. But neither of $D(n,n'), D(m,n)$ includes the other, since $m\notin D(n,n')$ and from the maximality of $D(n,n')$; and so $D(n,n')\cup D(m,n)$ is a closed disc with
the same boundary as $D(m,n')$. It is not equal to $D(m,n')$, from the maximality of $(D(n,n')$; so $\Sigma$ is a sphere
and these two are complementary discs. But then the union of $D(n,n'), D(m,n), D(m,n')$ equals $\Sigma$ and the claim holds.
This proves (1).

\bigskip
For each body $c$, either $X\cap \tilde{c}\ne \emptyset$, or $c$ is pinned by some pair of distinct vertices in $X$, since $X$ 
meets both wings of $c$ and their intersection is $\tilde{c}$. From (1), there is a set $Y$ of at most  $3|X|^2$ plates such that 
for every body $c$, either $X\cap \tilde{c}\ne \emptyset$, or $c$ is a subset of a plate in $Y$. This proves \ref{3planar2cell}.~\bbox

Now the complementary case:
\begin{thm}\label{3planar3cell}
For every surface $\Sigma$ without boundary,  there is a function $f$ such that, if $\Gamma$ satisfies {\bf (I1)}\LL {\bf (I6)},
 then for
all $k$, and for all choices of the set $\mac B$ of bodies such that every body has order three,
either there are $k$ pairwise disjoint $\mac B$-wings, or there is a set $X$ of nodes with $|X|\le f(k)$, and a set $Y$ of plates with $|Y|\le f(k)$,
such that for every body $c$, either $X\cap \tilde{c}\ne \emptyset$, or $c$ is a
subset of a plate in $Y$.
\end{thm}
\Proof Let $f_1$ be the function $f$ of \ref{hitwheels}, and define $f(k)=3f_1(k)^2$. We will show that $f$ satisfies the theorem.
Let $k\ge 0$ be an integer.
We assume that there do not exist $k$ pairwise disjoint $\mac B$-wings. Let $J$ be the graph drawn in $\Sigma$, with vertex set 
$N(\Gamma)$, and with each cell $c$ of order two replaced by an edge, and each cell of order three replaced by a triangle of edges, in the natural 
way. Let $U$ be the set of nodes of $\Gamma$ incident with bodies. Since there do not exist $k$ node-disjoint $\mac B$-wings in 
$\Gamma$, it follows that there do not exist $k$ vertices of $J$ in $U$ such that their wheel neighbourhoods in $J$ are pairwise disjoint.
By \ref{hitwheels}, there exists $X\subseteq N(\Gamma)$ with $|X|\le f_1(k)$, such that $X\cap W\ne \emptyset$
for every $\mac B$-wing $W$.

Let $x,x'\in X$ be distinct. We say a body $c$ is {\em pinned by} $\{x,x'\}$ if $x,x'\notin \tilde{c}$ and
there exist distinct $n,n'\in \tilde{c}$ such that $n$ sees $x$ and not $x'$, and $n'$ sees $x'$ and not $x$.
Let $\mathcal{C}$ be the set of all bodies pinned by $\{x,x'\}$, and let
$M$ be the union over all $c\in \mac C$ of the set of (two or three) $n\in \tilde{c}$ that can see exactly one of $x,x'$.
For each node $n\in M$ that can see $x$, there is a region incident with both $x,n$ and (therefore) not with $x'$;
choose a line segment $L(n,x)$ between $n,x$ with interior
within one such region. Choose these
all pairwise disjoint except for $x$.  For each $n\in M$ that can see $x'$, define $L(n',x')$ similarly. For each $c\in \mathcal{C}$, 
there exist $n,n'\in M$ such that $L(n,x)$ and $L(n,x')$ exist; we call the ordered pair $nn'$ a {\em tie}. 
For each tie $nn'$, let $L(n,n')$ be a 
line segment between $n,n'$ in the region incident with both $n,n'$ (which therefore is incident with neither of $x,x'$) choose all these line segments pairwise disjoint except for their ends.
\\
\\
(1) {\em If there exist ties $n_1n'$ and $n_2n'$ (with a common term $n'$), such that
$\ins(L(n_1,x)\cup L(n_1,n')\cup L(n_2,n')\cup L(n_2,x))$ includes $x'$, then there is a set $X_{xx'}$ with $X_{xx'}\le 1$ and a set $Y_{xx'}$ of at most three plates such that for every cell $c\in \mathcal{C}$, either $X_{xx'}\cap \tilde{c}\ne \emptyset$, or $c$ is a subset of a plate in $Y_{xx'}$.}
\\
\\
Suppose that there are two ties $n_1n'$ and $n_2n'$ with this property. 
Choose them with this disc ($D$ say) maximal. If $m,m'\in D$ for every tie $mm'$, then every body in $\mathcal{C}$
is inside $D$ except possible bodies incident with $n'$, and the claim holds.
So we may assume that 
there is a tie $mm'$ with one of $m,m'\notin D$. If one of $m,m'$ is in the interior of $D$, then the other is not
in $D$, and so $L(m,m')$
intersects one of $L(n_1,x), L(n_1,n'), L(n_2,n'), L(n_2,x)$, which is impossible (because we chose $L(m,m')$ disjoint from 
$L(n_1,n'), L(n_2,n')$, and it is inside a region not incident with either of $x,x'$, so it is disjoint from 
$L(n_1,x), L(n_2,x)$). So $m,m'$ are not in the interior of $D$. But $m'$ sees $x'$, which is in the interior of $D$; and since
each of $L(n_1,x), L(n_1,n'), L(n_2,n'), L(n_2,x)$ belongs to a region not incident with $x'$, it follows that $m'$ belongs to the
boundary of $D$, and hence is one of $n_1,n_2,n'$. But $m'\ne n_1,n_2$ since $m'$ sees $x'$ and $n_1,n_2$ do not; so $m'=n'$.
Thus $\ins(L(n_1,x)\cup L(n_1,n')\cup L(m,n')\cup L(m,x))$ and $\ins(L(m,x)\cup L(m,n')\cup L(n_2,n')\cup L(n_2,x))$ both bound
plates. Neither includes $D$, from the maximality of $D$, and so their union equals $\Sigma$ and the claim is true. 
This proves (1).
\\
\\
(2) {\em There is a set $X_{xx'}$ with $|X_{xx'}|\le 3$ and a set $Y_{xx'}$ of at most three plates such that for every cell $c\in \mathcal{C}$, either $X_{xx'}\cap \tilde{c}\ne \emptyset$, or $c$ is a subset of a plate in $Y_{xx'}$.}
\\
\\
For each tie $nn'$, $L(n,x)\cup L(n,n')\cup L(n',x')$ is a $\Gamma$-normal
line segment between $x,x'$ that passes through $n,n'$ and no other nodes except $x,x'$, and we denote it by $S_{nn'}$,
and call $S_{nn'}$ a {\em strut through $nn'$}.
We may assume that there are two disjoint ties, because otherwise there is a set of at most two nodes meeting all ties and the claim is true.
Hence there are two struts
$S_1,S_23$ with union a $\Gamma$-normal O-arc 
passing through at most six nodes, so its inside is defined. Choose $S_1,S_2$ with union an O-arc such that $\ins(S_1\cup S_2)$
is maximal. Let $D=\ins(S_1\cup S_2)$, and let $S_i$ be through 
$n_in_i'$ for $ i = 1,2$. If both $m,m'\in D$ for every tie $mm'$, then we are done as usual. So we assume without
loss of generality that there is a tie $mm'$ with $m\notin D$. As in the proof of (1), $m'$ is not in the interior of $D$.
Moreover, $m'\ne n_1,n_2$ since it does not see $x$. Thus either $m'\notin D$, or $m'$ is one of $n_1',n_2'$. From the symmetry we
may assume that $m'\ne n_2'$. Let $S$ be the strut through $mm'$. Hence $\ins(S\cup S_2)$ exists, and does not include 
$D$ from the maximality of $D$. If $m'\ne n_1'$ then similarly $\ins(S\cup S_1)$ exists, and the union of these three discs 
is $\Sigma$ and the claim holds. So we assume that $m'=n_1'$. The union of $D$ and $\ins(S\cup S_2)$ is a disc $D'$
with boundary the O-arc
$L(n_1,x)\cup L(n_1,n_1')\cup L(n_1',m)\cup L(m,x)$. Let $D''$ be the inside of this O-arc. By (1), $D'\ne D''$,
since $x'\in D'$. But then the three plates $D,D',D''$ satisfy the claim. This proves (2).

\bigskip

For each body $c$, either $X\cap \tilde{c}\ne \emptyset$, or $c$ is pinned by some pair of distinct vertices in $X$, since $X$
meets both wings of $c$ and their intersection is $\tilde{c}$. Let $X'$ be the union of the sets $X_{xx'}$ for all distinct $x,x'\in X$, and let $Y$ be similarly the union of the sets $Y_{x,x'}$. Thus $|X\cup X'|, |Y|\le 3|X|^2$. But 
for every body $c$, either $X\cap \tilde{c}\ne \emptyset$, or there $c$ is pinned by some pair of distinct vertices in $X$,
and then by (2), either $X_{xx'}$ intersects $\tilde{c}$ or $c$ is a subset of a plate in $Y_{xx'}$. 
Since $|X|\le f_1(k)$, this proves 
\ref{3planar2cell}.~\bbox

Then \ref{3planar} follows by partitioning $\mac B$ into two subsets, containing the cells of order two and those of order three, 
and applying \ref{3planar2cell} and \ref{3planar3cell} to them.


\section{Portraits in a surface without boundary}

We need a number of definitions.
A {\it separation} of a graph $G$ is a pair $(A,B)$ of subgraphs with $A \cup B = G$ and $E(A \cap B) =
\emptyset$, and its {\it order} is $|V(A \cap B)|$.  

We need the ``portrayals'' of~\cite{GM17}, but only in a surface without boundary, and only of graphs, not hypergraphs, so we will give them a different name, ``portraits''.
Let $G$ be a graph, and let $\Sigma$ be a connected surface without boundary.  A {\it portrait}
$\Pi = (\Gamma, \alpha,\gamma)$ of $G$ in $\Sigma$ consists of
\begin{itemize}
\item  a painting $\Gamma$ in $\Sigma$;
\item  a function $\alpha$ which assigns to each
cell $c$ of $\Gamma$ a subgraph $\alpha(c)$ of $G$; and 
\item  an injection $\gamma$ from 
$N(\Gamma)$ into $V(G)$
\end{itemize}
satisfying the axioms below.  For each $X \subseteq N(\Gamma)$ we denote $\{  \gamma (n) : n \in X \}$ by
$\gamma (X)$; for each cell $c_0 $ we denote $\bigcup ( \alpha (c) : c \in C(\Gamma) \setminus \{c_0\})$ by $\alpha (- c_0)$; and for $D\subseteq \Sigma$ we denote by $\alpha(D)$ the union of $\alpha(c)$ over all cells $c$ of $\Gamma$ with $c\subseteq D$.
The axioms are as follows:
\begin{itemize}
\item  $G = \alpha(\Sigma)$, and $E(\alpha (c) \cap \alpha (c')) = \emptyset$ for
all distinct cells $c, c'$; and
\item  $V(\alpha (c) \cap \alpha (-c))\subseteq \gamma ( \tilde{c}) \subseteq V (\alpha (c))$ for each cell $c$.
\end{itemize}
It is {\em matte} if $\alpha(c)$ is planar for every cell $c$; and {\em graphic} if $\Gamma$ has no cells of order three. 

If $\Pi=(\Gamma, \alpha,\gamma)$ is a portrait of $G$ in $\Sigma$, and $B$ is an induced subgraph of $G$, we can obtain a portrait 
$(\Gamma', \alpha',\gamma')$ of $B$ in $\Sigma$ as follows:
\begin{itemize}
\item For each cell $c$ of $\Gamma$, let $S=\{n\in \tilde{c}: \gamma(n)\in V(B)\}$. Choose a closed disc 
$D_c\subseteq \overline{c}$ such that $D_c\cap \tilde{c} = S$, taking $D_c=\overline{c}$ if $S=\tilde{c}$. 
\item Let $\Gamma'$ be the painting $(U,N)$ where $U$ is the union of $D_c$ over all cells $c$ of $\Gamma$, and 
$N=\{n\in N(\Gamma):\gamma(n)\in V(B)\}$. Thus every cell  $c'$ of $\Gamma'$ is included in a unique cell of $\Gamma$ which we
denote by $p(c')$. 
\item For each cell $c'$ of $\Gamma'$, define $\alpha'(c') = \alpha(p(c))\cap B$.
\item For each $n\in N(\Gamma')$, let $\gamma'(n) = \gamma(n)$.
\end{itemize}
It is straightforward to check that $\Pi'=(\Gamma', \alpha',\gamma')$ is indeed a portrait of $B$, and it is matte if 
$\Pi$ is matte. (Compare theorems (5.2), (5.4), (6.1), (6.2) and (6.4)
of~\cite{GM17}.) We write $\Pi' = \Pi[B]$, or $\Pi'=\Pi\setminus X$ where $X=V(G)\setminus V(B)$.

If $\Sigma$ is a surface, its boundary is the disjoint union of a number of O-arcs, which we call {\em cuffs}. 
If we paste a disc onto each cuff we obtain a surface $\hat{\Sigma}$ without boundary; this is called {\em capping} $\Sigma$.
If $F$ is an O-arc in a connected surface $\Sigma$ without boundary, and $F$ is not the boundary of a closed disc in $\Sigma$,
and we cut along $F$, we obtain a new surface $\Sigma_1$, possibly not connected, 
with boundary consisting of one or two O-arcs.  By capping $\Sigma_1$ 
we obtain a surface $\Sigma_2$ without boundary. Each of its components $\Sigma_3$ has
genus smaller than that of $\Sigma$. We say that such a surface $\Sigma_3$ (obtained from $\Sigma$ 
in this way by cutting and capping)
is {\em simpler}
than $\Sigma$. Any portrait in $\Sigma_3$ can be converted to a portrait in $\Sigma$ in the natural way.

In this section we prove the following:
\begin{thm}\label{3planar2}
For every connected surface $\Sigma$ without boundary, there is a  function $g$ such that for all integers $k\ge 0$, 
if $G$ admits a matte portrait in $\Sigma$, then either:
\begin{itemize}
\item $G$ has K-number at least $k$; or
\item there exists $Z\subseteq V(G)$ with $|Z|\le g(k)$ such that $G\setminus Z$ can be drawn in $\Sigma$; or
\item there is a separation $(A,B)$ of $G$ of order at most $g(k)$ such that both $A,B$ are not planar.
\end{itemize}
\end{thm}

We will derive it from \ref{3planar2cell} and \ref{3planar3cell}. If  $(\Gamma, \alpha, \gamma)$ is a portrait of $G$ in $\Sigma$,
we say a cell $c$ of $\Gamma$ is {\em flat} if $\alpha(c)$ can be drawn in a closed disc such that the vertices
$\gamma(n)\;(n\in \tilde{c})$ are drawn in the boundary of the disc.
It is convenient to add two more conditions 
about a portrait $(\Gamma, \alpha, \gamma)$ to the list 
{\bf (I1)}--{\bf (I6)}, namely:
\begin{itemize}
\item[{\bf (I7)}] Let $c$ be a cell, and let
$v_1 ,v_2 \in \gamma (\tilde{c})$.  Then there is a path
in $\alpha(c)$ from $v_1$ to $v_2$ with
no  other vertex in $\gamma(\tilde{c})$.
\item[{\bf (I8)}] For every cell $c$ of order three, $\alpha(c)$ is not flat.
\end{itemize}

\begin{thm}\label{setup}
For every connected surface $\Sigma$ without boundary, 
if $G$ admits a matte portrait $\Pi_0$ in $\Sigma$, then either:
\begin{itemize}
\item there exists $Z\subseteq V(G)$ with $|Z|\le 6$ and a connected surface $\Sigma'$ without boundary that is simpler than 
$\Sigma$, such that $G\setminus Z$ admits a matte portrait in $\Sigma'$, graphic if $\Pi_0$ 
is graphic; or
\item there is a separation $(A,B)$ of $G$ of order at most $6$ such that both $A,B$ are not planar; or
\item there exists $Z\subseteq V(G)$ with $|Z|\le 6$ such that $G\setminus Z$ is planar; or 
\item there is a matte portrait $\Pi$ in $\Sigma$, satisfying {\bf (I1)}--{\bf (I8)}, and graphic if $\Pi_0$ is graphic.
\end{itemize}
\end{thm}
\Proof
We proceed by induction on $|V(G)|$. 
\\
\\
(1) {\em We may assume that $G$ is 2-connected.}
\\
\\
Suppose that $(A,B)$ is a separation of order at most one, with $V(A),V(B)\ne V(G)$. We may assume that $A$ is planar, because otherwise
the second outcome holds.
Let $\Pi_0=(\Gamma, \alpha, \gamma)$ be a matte portrait
of $G$ in $\Sigma$.
Thus $\Pi_0[B]=(\Gamma', \alpha', \gamma')$
is a matte portrait of $B$ in
$\Sigma$, graphic if $\Pi_0$ is graphic. From the inductive hypothesis, either:
\begin{itemize}
\item there exists $Z\subseteq V(B)$ with $|Z|\le 6$ and a connected surface $\Sigma'$ without boundary that is simpler than 
$\Sigma$, such that $B\setminus Z$ admits a matte portrait in $\Sigma'$, graphic if $\Pi_0$
is graphic (but then the same is true for $G\setminus Z$, since $A\setminus Z$ is planar and contains at most one vertex of $B\setminus Z$);
\item there is a separation $(A',B')$ of $B$ of order at most $6$ such that both $A',B'$ are not planar (but then
one of $(A\cup A', B')$, $(A', A\cup B')$ is a separation of $G$ with the same order as $(A',B')$, and both its terms are nonplanar); or
\item there exists $Z\subseteq V(G)$ with $|Z|\le 6$ such that $B\setminus Z$ is planar (but then $G\setminus Z$ is planar, since
$A$ is planar and contains at most one vertex of $B$); or
\item there is a matte portrait $\Pi$ of $B$ in $\Sigma$, satisfying {\bf (I1)}--{\bf (I8)}, and graphic if $\Pi_0$ is graphic (but then the same is true for $G$ since $A$ is planar and contains at most one vertex of $B$).
\end{itemize}
In each case the theorem holds. This proves (1).
\\
\\
(2) {\em We may assume that, for every matte portrait $\Pi=(\Gamma, \alpha, \gamma)$
of $G$ in $\Sigma$, graphic if $\Pi_0$ is graphic, and 
 every $\Gamma$-normal O-arc $F$ in $\Sigma$ that contains at most six nodes of $\Gamma$, there is a 
closed disc $\Delta$ in $\Sigma$ with boundary $F$ such that $\alpha(D)$ is planar.}
\\
\\
Suppose that $F$ is a $\Gamma$-normal O-arc that contains at most six nodes of $\Gamma$, and $F$ is not the boundary of such a disc.
Let $N=F\cap N(\Gamma)$, and $Z=\gamma(N)$. Let $\Sigma'$ be obtained from $\Sigma$ by cutting along $F$ and capping. Let
$\Sigma'$ have components $\Sigma_1\LL \Sigma_t$ say (thus $t\le 2$). For $1\le i\le t$, let $G_i=\alpha(\Sigma_i)$.
If more than one of $G_1\LL G_t$ are not planar, then the second outcome of the theorem holds. So we may assume that $G_2\LL G_t$ 
are planar.
Now $\Pi\setminus Z$
is a matte portrait of $G\setminus Z$ in $\Sigma$, graphic if $\Pi_0$ is graphic, and $F$ is disjoint
from the closures of all its cells; and hence there is a matte portrait 
of $G_1\setminus Z$ in $\Sigma_1$, graphic if $\Pi_0$ is graphic. Since $G_2\LL G_t$
are planar, it follows that there is a matte portrait 
$(\Gamma', \alpha', \gamma')$ of $G\setminus Z$ in $\Sigma_1$ (by redrawing the planar graphs $G_2\LL G_t$ in $\Sigma_1$ appropriately).

If $\Sigma_1$ is simpler than $\Sigma$, 
the first outcome holds. So
we may assume that $\Sigma_1$ is not simpler than $\Sigma$. Hence $t=2$, and $\Sigma_2\subseteq \Sigma$ is a disc with boundary $F$.
Since $G_2$ is planar, this proves (2).

\bigskip

For $i=0,1,2,3$ let $n_i$ be the number of cells of $\Gamma$ with order $i$. 
Let us
choose $\Gamma, \alpha, \gamma$, with $n_3$ minimum; subject to that, with $n_2$ minimum; 
subject to that, with $n_1$ minimum; and subject to that, with $n_0$ minimum. We call this the ``optimality''
of the portrait. We may assume that $G$ is not planar.  
\\
\\
(3) {\em The following hold:
\begin{itemize}
\item $\alpha(c)$ is nonnull for every cell $c$ (because otherwise $c$ could be removed from $\Gamma$, contrary to optimality);
\item $\alpha(c)\ne G$ for each cell $c$ (since $\alpha(c)$ is planar and $G$ is not);
\item $n_0=0$ (because $\alpha(c)$ is nonnull and not equal to $G$, and $G$ is connected);
\item $n_1=0$ (because if $\tilde{c}=\{n\}$ say, then $\alpha(c)$ has only one vertex since $G$ is 2-connected by (1), and 
so $c$ can be removed from $\Gamma$, contrary to optimality).
\end{itemize} 
}

If $F$ is a $\Gamma$-normal O-arc passing through at most six nodes, let $\ins(F)$ be a disc as in (2). Such a disc is a {\em plate}.
\\
\\
(4) {\em Let $F$ be a $\Gamma$-normal O-arc passing through at most three nodes.
\begin{itemize}
\item If $|F\cap N(\Gamma)|\le 1$ then $\ins(F)\cap U(\Gamma) = F\cap N(\Gamma)$.
\item If $|F\cap N(\Gamma)|=2$ then $\ins(F)$ includes at most one cell, and if so then that cell $c$ satisfies $\tilde{c} =F\cap N(\Gamma)$.
\item If $|F\cap N(\Gamma)|=3$ then either $\ins(F)$ includes no cells of order three, or it includes exactly one cell, and that cell 
satisfies $\tilde{c} =F\cap N(\Gamma)$.
\end{itemize}}
\noindent
The proofs are clear from optimality.
\\
\\
(5) {\em Every region is an open disc, and the closure of every region
is a closed disc.}
\\
\\
(Compare theorems (8.1) and (8.2) of~\cite{GM17}.) For every O-arc $F$ in $r$, it bounds a plate with interior disjoint from 
$U(\Gamma)$ (by (4)), and hence included in $r$; and so $r$ is an open disc. Moreover, for every $\Gamma$-normal O-arc $F$
with $|F\cap N(\Gamma)|=1$, again by (4) $F$ bounds a disc included in $r\cup F$; and so $\overline{r}$ is a closed disc.
This proves (5).
\\
\\
(6) {\em We may assume that $\Gamma$ satisfies {\bf (I1)}--{\bf (I6)}.}
\\
\\
From (2)--(5), it follows that $\Gamma$ satisfies {\bf (I1)}\LL {\bf (I3)} and {\bf (I6)}. For {\bf (I4)},
let $c$ be a cell of order two, and let $r,r'$ be the regions incident with $c$; let $\tilde{c}=\{n_1.n_2\}$, and suppose that
$n_3\ne n_1,n_2$ is incident with $r_1,r_2$. For $1\le i<j\le 3$, choose a line segment $L(i,j)$ with interior in $r$, 
between $n_i,n_j$, disjoint except for their ends; and choose $L'(i,j)\subseteq r'$ similarly. Then for $1\le <j\le 3$, $L(i,j)\cup L'(i,j)$ is an O-arc $F_{ij}$ say; let $D_{ij}$ be its inside. From the second statement of (4), none of $D_{12}, D_{13}, D_{23}$
includes either of the other two, and so their union, togther with $r,r'$, equals $\Sigma$. Consequently 
$G\setminus \{\gamma(n_1),\gamma(n_2), \gamma(n_3)\}$ is planar and the theorem holds. Thus, we may assume that
$\Gamma$ satisfies {\bf (I4)}.

For  {\bf (I5)}, $c$ be a cell of order three, with $\tilde{c}=\{n_1,n_2,n_3\}$. Let $n_4\ne n_1,n_2,n_3$, and suppose that
for $i = 1,2,3$, there is a region $r_i$ incident with both $n_i, n_4$.
For $1\le i\le 3$, let $L(i,4)$ be a line segment between $n_i,n_4$ with interior in $r_i$, and with interior disjoint from $L(j,4)$
for $j\in \{1,2,3\}\setminus \{i\}$ (this last is automatic unless $r_i=r_j$). For $1\le i<j\le 3$, let $L(i,j)$ be a line segment
between $n_i,n_j$ and with interior in the region incident with $n_i,c$ and $n_j$, with interior disjoint from $L_{i,4}$
and $L(j,4)$. For $1\le i<j\le 3$, $L(i,4)\cup L(j,4)\cup L(i,j)$ is an O-arc $F_{ij}$; let $D_{ij}$ be its inside. Let $F_4$
be the O-arc $L(1,2)\cup L(1,3)\cup L(2,3)$, and let $D_4$ be its inside. None of $D_1,D_2,D_3$ contains $n_i$, by the third 
statement of (4). So their union, together with $D_4$, equals $\Sigma$, and the theorem holds.  This proves (6).
\\
\\
(7) {\em Let $c$ be a cell, and let
$v_1 ,v_2 \in \gamma (\tilde{c})$.  Then there is a path
in $\alpha(c)$ from $v_1$ to $v_2$ with
no  other vertex in $\gamma(\tilde{c})$. Consequently, $\Pi$ satisfies {\bf (I7)}.}
\\
\\
This is proved like theorem (9.1) of~\cite{GM17}.
\\
\\
(8) {\em If $c$ is a cell of order three, then $c$ is not flat. Consequently, $\Pi$ satisfies {\bf (I8)}.}
\\
\\
If $\alpha(c)$ is flat, then we can draw $\alpha(c)$ in the disc $\overline{c}$ with $\gamma(n)$
mapped to $n$ for each $n\in \tilde{c}$; and then, by thickening every edge of this drawing to a cell of order two incident with the
same two nodes, we obtain a matte portrait of $G$ in $\Sigma$ contrary to optimality. This proves (8).

\bigskip

From (6), (7), (8), this completes the proof of \ref{setup}.~\bbox

\begin{thm}\label{disjtframes}
Let $(\Gamma, \alpha, \gamma)$ be a matte portrait of $G$ in $\Sigma$, satisfying
{\bf (I1)}--{\bf (I8)}; and let $\mathcal{B}$ be a subset of the set of all cells of $\Gamma$, such that no cell in $\mac B$ is flat. 
For all $k$, if there are $k$ pairwise disjoint $\mathcal{B}$-wings in $\Gamma$, then there are $k$ vertex-disjoint K-graphs in $G$.
\end{thm}
\Proof
Let $c_0\in \mac B$, and let $W$ be a wing of $c_0$. We define a {\em frame } of $W$ as follows.

Assume first that $c_0$ has order two; so $W$ is the closure of a region $r$ incident with $c_0$. For each cell $c$ incident with $r$ with $c\ne c_0$,
there are two nodes $m_c,n_c$ incident with $c$ and with $r$; let $P_c$ be a path between $\gamma(m_c), \gamma(n_c)$ as in {\bf (I7)}.
The union of all the paths $P_c$ is a path with ends $\gamma(n_1), \gamma(n_2)$, where $\tilde{c_0}=\{n_1,n_2\}$; and we call this path
a {\em frame} for the wing $W=\overline{r}$. Note that since $c_0$ is not flat, the union of $\alpha(c_0)$ and the frame for $r$ is nonplanar, and so
contains a K-graph.

Now we assume that $c_0$ has order three; so for some $n_1\in \tilde{c}$,  $W$ is the union of the closures of all regions
incident with $n_1$. For each region $r$ incident with $n_1$, and for every cell $c$ incident with $r$, let $m_c,n_c$ be the
two nodes incident with $c$ and with $r$, and let $P_c$ be a path between $\gamma(m_c), \gamma(n_c)$ as in (8). The union of all
the paths $P_c$ (over all regions $r$ incident with $n_1$, and all cells $c$ incident with $r$) is a connected subgraph, that contains
$\gamma(n)$ for all $n\in \tilde{c}$, and we call it a {\em frame} for the wing $B$. Moreover, let $\tilde{c}=\{n_1,n_2,n_3\}$; then
every two of $n_1,n_2,n_3$ can be joined by a path in the frame that does not pass through the third member of $\tilde{c}$.
Consequently, since $c_0$ is not flat,  the union of $\alpha(c)$ and the frame is nonplanar, and so contains a K-graph.

Suppose that $W_1\LL W_k$ are $\mathcal{B}$-wings, pairwise disjoint. For $1\le i\le k$, let $W_i$ be a wing of a cell 
$c_i\in \mac B$, and let $F_i$ be a frame for $W_i$. 
We claim that $\alpha(c_i)\cup F_i\;(1\le i\le k)$ are pairwise vertex-disjoint. Suppose some vertex $v$
belongs to $\alpha(c_1)\cup F_1$ and to $\alpha(c_2)\cup F_2$ say. Since $c_1\ne c_2$ (because $B_1,B_2$ are disjoint, and each contains
two nodes incident with $c$), it follows that $v\in V(F_1\cap F_2)$. If $v=\gamma(n)$ for some node $n$, then $n\in B_1\cap B_2$,
which is impossible since $B_1\cap B_2=\emptyset$. So $v\in V(\alpha(c))\setminus \tilde{c}$, for some cell $c$, and there is no other
cell $c'$ with $v\in V(\alpha(c'))$. Since $v\in V(F_1)$, $c$ is incident with a region included in $B_1$, and similarly, with a
incident with a region included in $B_2$. But no cell is incident with a region in $B_1$ and with a region in $B_2$ since $B_1,B_2$
are disjoint. Consequently $G$ has $k$ pairwise vertex-disjoint K-graphs and the theorem holds. This proves \ref{disjtframes}.~\bbox

We denote the genus of a surface $\Sigma$ by $\vare(\Sigma)$. The next step is:
\begin{thm}\label{3planar22}
For every connected surface $\Sigma$ without boundary, there is a  function $g$ such that for all integers $k\ge 0$,
if $G$ admits a graphic matte portrait in $\Sigma$, then either:
\begin{itemize}
\item $G$ has K-number at least $k$; or
\item there exists $Z\subseteq V(G)$ with $|Z|\le g(k)$ such that $G\setminus Z$ can be drawn in $\Sigma$; or
\item there is a separation $(A,B)$ of $G$ of order at most $6\vare(\Sigma)+6$ such that both $A,B$ are not planar.
\end{itemize}
\end{thm}
\Proof Let $f_1$ be the function $f$ of \ref{3planar2cell}.
We proceed by induction on $\vare(\Sigma)$; so we assume the result holds for all simpler connected surfaces without boundary.
Choose a function $g'$ that satisfies the theorem with $g$ replaced by $g'$, for all $\Sigma'$ simpler than $\Sigma$.
Let $g(k)=7f_1(k)+g'(k)+6$. We will show that $g$ satisfies the theorem.

By \ref{setup}, either:
\begin{itemize}
\item there exists $Z\subseteq V(G)$ with $|Z|\le 6$ and a connected surface $\Sigma'$ without boundary that is simpler than 
$\Sigma$, such that $G\setminus Z$ admits a graphic matte portrait in $\Sigma'$; or
\item there exists $Z\subseteq V(G)$ with $|Z|\le 6$ such that $G\setminus Z$ is planar; or
\item there is a separation $(A,B)$ of $G$ of order at most $6$ such that both $A,B$ are not planar; or
\item there is a graphic matte portrait of $G$ in $\Sigma$, satisfying {\bf (I1)}--{\bf (I8)}.
\end{itemize}
In the first case, the inductive hypothesis, applied to $G\setminus Z$ and $\Sigma'$, tells us that 
either:
\begin{itemize}
\item $G\setminus Z$ has K-number at least $k$ (and then so does $G$); or
\item there exists $Z'\subseteq V(G\setminus Z)$ with $|Z'|\le g'(k)$ such that $G\setminus (Z\cup Z')$ can be drawn in $\Sigma'$
(and then the second outcome of the theorem holds, since $g(k)\ge g'(k)+6$, and every portrait in $\Sigma'$ can be converted to a 
portrait in $\Sigma$); or
\item there is a separation $(A',B')$ of $G\setminus Z$ of order at most $6\vare(\Sigma')+6$ such that both $A',B'$ are not planar
(but then there is a separation $(A,B)$ of $G$ with $Z\subseteq V(A\cap B)$, and $A\setminus Z=A'$, and $B\setminus Z=B'$, and so 
$(A,B)$ satisfies the third outcome of the theorem since its order is at most $6\vare(\Sigma')+12\le 6\vare(\Sigma)+6$).
\end{itemize}
In the second case, the second outcome of the theorem holds, and in the third case, the third outcome of the theorem holds. 
Consequently we may assume that the fourth case holds.

Let $(\Gamma, \alpha, \gamma)$ be a graphic matte portrait of $G$ in $\Sigma$, satisfying {\bf (I1)}--{\bf (I8)}.
Let $\mac B$ be the set of all cells of $\Gamma$ that are not flat. 
By \ref{disjtframes}, we may assume that there do not exist $k$ pairwise disjoint $\mac B$-wings. Hence, by \ref{3planar2cell},
there is a set $X$ of nodes with $|X|\le f_1(k)$, and a set $Y$ of plates with $|Y|\le f_1(k)$,
such that for every body $c$, either $X\cap \tilde{c}\ne \emptyset$, or $c$ is a
subset of a plate in $Y$. 
For each plate $D\in Y$, $\alpha(D)$
is planar. Let $S$ be the set of nodes of $\Gamma$ that belong to the boundary of a plate in $Y$; so $|S|\le 6|Y|\le 6f_1(k)$. 
Let $Z=\{\gamma(n):n\in S\}$. Let $A$ be the union of $\alpha(D)$ over all $D\in Y$; then each component of $A\setminus Z$ is planar, 
because it is a subgraph of some $\alpha(D)$. Consequently $A\setminus Z$ is planar. 
Let $B'=B\setminus (\gamma(X)\cup Z)$. Then $\Pi[B']=(\Gamma', \alpha', \gamma')$ is a matte portrait in $\Sigma$.
But every cell of $\Gamma'$ is flat, since
cells of order at most one are always flat in a matte portrait, and cells of $\Gamma'$ of order two are flat, because every 
such cell $c$ is a cell 
of $\Gamma$ with $\tilde{c}\cap X=\emptyset$. It follows that $B\setminus (\gamma(X)\cup Z)$ can be drawn in $\Sigma$,
and therefore so can $G\setminus (\gamma(X)\cup Z)$. This proves \ref{3planar22}.~\bbox

Armed with \ref{3planar22} we can prove \ref{3planar2}, 
which we restate:
\begin{thm}\label{3planar2again}
For every connected surface $\Sigma$ without boundary, there is a  function $g$ such that for all integers $k\ge 0$,
if $G$ admits a matte portrait in $\Sigma$, then either:
\begin{itemize}
\item $G$ has K-number at least $k$; or
\item there exists $Z\subseteq V(G)$ with $|Z|\le g(k)$ such that $G\setminus Z$ can be drawn in $\Sigma$; or
\item there is a separation $(A,B)$ of $G$ of order at most $g(k)$ such that both $A,B$ are not planar.
\end{itemize}
\end{thm}
\Proof
The proof is very much like that of \ref{3planar22}. Let $f_2$ be the function $f$ of \ref{3planar3cell}, and let $g_1$ be the function of \ref{3planar22}.
We proceed by induction on $\vare(\Sigma)$; so we assume the result holds for all simpler connected surfaces without boundary.
Choose a function $g'$ that satisfies the theorem with $g$ replaced by $g'$, for all $\Sigma'$ simpler than $\Sigma$.
Let $g(k)=7f_2(k)+g'(k)+6+g_1(k)+6\vare(\Sigma)$. We will show that $g$ satisfies the theorem.

As in the proof of \ref{3planar22}, by applying \ref{setup}, we may assume that 
there is a matte portrait $(\Gamma, \alpha, \gamma)$ of $G$, satisfying {\bf (I1)}--{\bf (I8)}.
Let $\mac B$ be the set of all cells of $\Gamma$ of order three that are not flat.
By \ref{disjtframes}, we may assume that there do not exist $k$ pairwise disjoint $\mac B$-wings. Hence, by \ref{3planar3cell},
there is a set $X$ of nodes with $|X|\le f_2(k)$, and a set $Y$ of plates with $|Y|\le f_2(k)$,
such that for every body $c$, either $X\cap \tilde{c}\ne \emptyset$, or $c$ is a
subset of a plate in $Y$.
Thus $\alpha(D)$
is planar for each $D\in Y$. Let $S$ be the set of nodes of $\Gamma$ that belong to the boundary of a plate in $Y$; so $|S|\le 6|Y|\le 6f_1(k)$.
Let $Z=\{\gamma(n):n\in S\}$. Let $A$ be the union of $\alpha(D)$ over all $D\in Y$; then
$A\setminus Z$ is planar.
Let $B'=B\setminus (\gamma(X)\cup Z)$. Then $\Pi[B']=(\Gamma', \alpha', \gamma')$ is a matte portrait in $\Sigma$.
Suppose that $c$ is a cell of order three in $\Pi[B']$. From the definition of $\Pi[B']$ it follows that $c$ is a cell of $\Gamma$,
and $\alpha(c) = \alpha'(c)$. But by {\bf (I8)}, $c$ is not flat in $\Pi$, and so $c\in \mathcal{B}$, and so 
either $X\cap \tilde{c}\ne \emptyset$ or $c$ is a subset of a plate in $Y$. In either case this contradicts that $c$ is a cell
of order three in $\Pi[B']$. So $\Pi[B']$ is graphic. By \ref{3planar22} applied to $\Pi[B']$, either:
\begin{itemize}
\item $B'$ has K-number at least $k$ (but then so does $G$); or 
\item there exists $Z'\subseteq V(B)$ with $|Z'|\le g_1(k)$ such that $B\setminus Z'$ can be drawn in $\Sigma$ (but then
so can $G\setminus (\gamma(X)\cup Z\cup Z')$, since $A$ is planar); or
\item there is a separation $(C',D')$ of $B'$ of order at most $6\vare(\Sigma)$ such that both $C', D'$ are not planar;
but then there is a separation $(C,D)$ of $G$ of order at most $6\vare(\Sigma)+ 7f_2(k)\le g(k)$ such that $C,D$ are not planar.
\end{itemize}
This proves \ref{3planar2again}.~\bbox

\section{Societies and crosses}

A {\em society} is a pair $(G,\Omega)$, where $G$ is a graph and $\Omega$ is a cyclic permutation of a subset 
(denoted by $\overline{\Omega}$) of $V(G)$; and a {\em cross} in a society $(G,\Omega)$ is a subgraph $P_1\cup P_2$
with two components $P_1,P_2$, each with both ends in $\overline{\Omega}$, and such that the ends alternate in $\Omega$;
that is, if $P_i$ has ends $s_i,t_i$ for $i = 1,2$, then $s_1,s_2,t_1,t_2$ occur in $\Omega$ in that order or its reverse.

Let us take a portrait $(\Gamma, \alpha, \gamma)$ of a graph $G$ in a sphere $\Sigma$, and choose a $\Gamma$-normal 
O-arc $F$ that bounds a closed disc including $U(\Gamma)$.
If we enumerate $N(\Gamma)\cap F$ in cyclic order, this induces a cyclic permutation of the set $\gamma(N(\Gamma)\cap F)$, 
say $\Omega$. Let us say such a portrait (with such a O-arc $F$) is a 
a {\em rural portrait} of $(G,\Omega)$. We need the following:
\begin{thm}\label{nocross}
A society has no cross if and only if it admits a rural portrait.
\end{thm}
For ``internally 4-connected'' societies, this is proved in section 2 of~\cite{GM9}, and for general societies it follows from 
theorems (11.6), (11.9) and (11.10) of the same paper.

A surface $\Sigma$ is an {\em annulus} if it can be obtained from a sphere by deleting the interiors of two disjoint closed discs.
Let $F,F'$
be disjoint O-arcs, with union the boundary of an annulus $\Sigma$.
An O-arc $R\subseteq \Sigma$ {\em separates} $F,F'$ if it bounds two closed discs
in $\hat{\Sigma}$, one including $F$ and the other including $F'$. 

Let $\Gamma$ be a painting of a graph $G$ in a surface $\Sigma$. If $c$ is a cell, each component of $\overline{c}\setminus \tilde{c}$ is the interior of a line segment with both ends in $\tilde{c}$, and we call that a {\em side} of $c$. 
An O-arc $F$ is {\em $\Gamma$-borderline} if $F$ is a subset of the union of the boundaries of the cells of $\Gamma$, 
and includes at most one side of each cell.
Let $\Theta$ be a cuff of $\Sigma$. An O-arc $F$ {\em surrounds} $\Theta$ if $F\cap \Theta=\emptyset$ and there is an annulus 
$\Sigma'\subseteq \Sigma$ with boundary $F\cup \Theta$.
For each node $n$ of $\Gamma$, let $d(n, \Theta)=1$ if $n\in \Theta$, and otherwise let $d(n, \Theta)$ be the minimum of 
$|L\cap N(\Gamma)|$ over all 
$\Gamma$-normal line segments $L$ between $n$ and some node in $\Theta$.
If $k\ge 2$ is an integer, a {\em $k$-ring around $\Theta$} means a $\Gamma$-borderline O-arc $R$ that surrounds $\Theta$, such that 
$d(n, \Theta) = k$ for every node in $R$.
A {\em $k$-nest around $\Theta$} is a sequence $(R_1\LL R_k)$ of O-arcs, pairwise disjoint, such that 
$R_1=\Theta$ and $R_i$ is an $i$-ring around $\Theta$ for $2\le i\le k$. We need:

\begin{thm}\label{getnest}
Let $\Gamma$ be a painting in a surface $\Sigma$, let $\Theta$ be a cuff of $\Sigma$, and let $F$ be a $\Gamma$-normal O-arc,
such that $d(n,\Theta)\ge k$ for every node $n\in F$. Then there is a $k$-nest around $\Theta$.
\end{thm}

The proof is clear and we omit it.
Next we need the following lemma:

\begin{thm}\label{crosstoK}
Let $(A,B)$ be a separation of a graph $G$. 
Let $B$ be drawn in a closed disc $D$, and let the vertices of $B$ drawn in $\bd (D)$ be $v_1\LL v_T$ in order. We will 
identify $B$ with the drawing of $B$, for simplicity. Suppose that
\begin{itemize}
\item all edges of $B$ are drawn within the interior of $D$;
\item $V(A\cap B)=\{v_t,v_{t+1}\LL v_T, v_1\}$;
\item there are paths $P_1,P_2,P_3$ of $G$, pairwise vertex-disjoint, where $P_i$ has ends $v_i$ and $v_{t+1-i}$ for $i = 1,2,3$;
$P_1$ is a path of $A\cup B$, and $P_2,P_3$ are paths of $B$;
\item every vertex in $A\cap B$ belongs to $P_1$; 
\item there exist $3\le a<b<c<d\le t-2$ such that there are two disjoint paths $Q,R$ of $G$ with ends $v_a,v_c$ and $v_b,v_d$ 
respectively.
\end{itemize}
Then $G$ is not planar.
\end{thm}
\Proof We proceed by induction on $|V(G)|+|E(G)|$. We may therefore assume that $G$ is the union of $P_1,P_2,P_3,Q,R$,
Also, we may assume that no edge with an end different from $v_1\LL v_t$ belongs to one of $P_1,P_2,P_3$ and to one of $Q,R$, because we could 
contract any such edge.
We may also assume that no internal vertex of $Q$ or $R$ is one of $v_3\LL v_{t-2}$, because if say $v_i$
is an internal vertex of $Q$, then we can replace $Q$ by one of its maximal subpaths with one end $v_i$ and still satisfy the 
hypotheses. Consequently, every edge belongs to exactly one of the five paths, and 
each of $v_1\LL v_t$ is an end of one of these five paths. 
Since each of $Q,R$ contains at least two vertices of each of $P_1,P_2,P_3$; it follows that $P_1,P_2,P_3$ all have length at least three, and $Q,R$ have length at least five. 
Hence for every edge $uv$ of $B$,
we may assume that either $u,v$ belong to distinct paths among $P_1,P_2,P_3$, or one is in $V(Q)$ and the other in $V(R)$,
because otherwise we could contract the edge $uv$ and apply the inductive hypothesis.
In particular, $v_3,v_{t-2}$ both belong to $V(Q\cup R)$, and are not internal vertices of these paths;
so $t=8, a=3,b=4, c=5$ and $d=6$. Similarly, $v_4,v_5$ both belong to $V(P_3)$. There is a disc in $D$ bounded by  
$P_3$ and part of the boundary of $D$, that contains $P_1$ and $P_2$; and it follows that every edge of $B$ is drawn within this disc, since every edge of $Q\cup R$ joins two of 
$P_1,P_2,P_3$. Thus $P_3$ has length three, since otherwise there is a vertex $v$ of $P_3$ different from $v_3,v_4,v_5,v_6$,
and we could add it to the sequence $v_1\LL v_8$ in the appropriate position to contradict that $t=8$. 
Now $Q$ has ends $v_3,v_5$; let $q$ be the 
neighbour of $v_5$ in this path. Similarly let $r$ be the neighbour of $v_4$ in $R$. It follows that $q,r\in V(P_2)$ 
(because every edge of $Q\cup R$ has ends in different paths among $P_1,P_2,P_3$), and for the same reason, no vertices or 
edges of $G$ are drawn in the interior of the disc of $D$ bounded by the cycle $C$, where $C$ is the union of the 
path $q\DD v_4\DD v_5\DD r$
and the subpath of $P_2$ between $r,q$.  Let $Q'$ consist of the union of the subpath of $Q$ between $q,v_3$ and the edge $v_3v_4$;
and let $R'$ similarly consist of the union of the subpath of $R$ between $r,v_6$ and the edge $v_5v_6$. 
Now, $G\setminus V(C)$ is connected,
since $G$ is the union of $P_1,P_2,P_3,Q,R$, and every edge of $Q\cup R$ has an end in $V(P_3)$ except for their first and
last edges. Hence
if $G$ has a drawing in a sphere then $C$ bounds a region in that drawing;
but that is impossible because $Q',R'$ are disjoint paths, and their ends are in $V(C)$ in alternating order.
So $G$ is not planar. This proves \ref{crosstoK}.~\bbox

\section{Tangles}

A {\it tangle} of {\it order} $\theta$ in
a graph $G$, where $\theta \ge 1$ is an
integer, is a set $\mathcal T$ of separations of $G$, all of order $< \theta$, such that:
\begin{enumerate}
\item
for every separation $(A, B)$ of order $< \theta$, one of $(A,B), (B,A)$ belongs to $\mathcal T$;
\item
if $(A_1, B_1), (A_2, B_2), (A_3, B_3) \in \mathcal T$ then $A_1 \cup A_2 \cup A_3\not= G$; and
\item
if $(A, B)\in \mathcal T$ then $V(A)\not= V(G)$.
\end{enumerate}
We define $ord (\mathcal T) = \theta$.

If $Z \subseteq V(G)$, we denote the graph obtained by deleting $Z$ by $G\setminus Z$.  If $\mathcal T$ is a
tangle in $G$ of order $\theta$ and $Z \subseteq V(G)$ with $|Z| < \theta$, we denote by
$\mathcal T\setminus Z$ the set of all separations $(A',  B')$ of $G\setminus Z$ of order $< \theta - |Z|$
such that there exists $(A, B)\in \mathcal T$ with $Z \subseteq V(A \cap B), A \setminus Z = A'$ and $B\setminus Z
= B'$.  It is shown in theorem (8.5) of \cite{GM10} that $\mathcal T\setminus Z$ is a
tangle in $G\setminus Z$ of order $\theta - |Z|$.
A tangle $\mac T$ in a graph $G$ is {\em matted} if $A$ is planar for each $(A,B)\in \mac T$.
\section{General portrayals}

Now we have to work with paintings and portrayals in surfaces with boundary, so we need a number of new definitions. 
Let $\Sigma$ be a surface (with boundary). A {\em painting} $\Gamma$ in $\Sigma$ is a pair $(U,N)$, satisfying the same 
conditions as in a surface without boundary, and in addition:
\begin{itemize}
\item $\bd (\Sigma) \subseteq U$
\item for each cell $c$, if $c \cap \bd ( \Sigma ) \neq \emptyset$
then $| \tilde{c} | = 2$ and
$\overline{c} \cap \bd (\Sigma )$ is a line segment
with ends the two members of $\tilde{c}$.
\end{itemize}
We use the same terminology and notation as before for paintings. Cells $c$ with $c \cap \bd ( \Sigma ) \neq \emptyset$
are {\em border} cells, and the others are {\em internal} cells. Similarly, nodes in $\bd(\Sigma)$ are {\em border} nodes, and 
the others are {\em internal} nodes.

Portrayals were defined in~\cite{GM17}; and the bulk of that paper was showing that if $G$ admits a portrayal in a surface 
$\Sigma$, then we can delete a few vertices such that what remains
admits a portrayal, in $\Sigma$ or a simpler surface, with several additional desirable properties. We only need portrayals 
in this form, so we will skip the general definitions and go straight to these prettified portrayals, which we call ``clean''.

A {\em clean portrayal} is a 5-tuple $(\Sigma, \Gamma, \alpha,\beta,\gamma)$, where:
\begin{itemize}
\item $\Sigma$ is a non-null connected surface with (possibly null) boundary;
\item $\Gamma$ is a painting $\Gamma$ in $\Sigma$, in which every cell has order two or three;
\item$ \alpha$ is a function which assigns to each cell $c$ of $\Gamma$ a subgraph $\alpha(c)$ of $G$:
\item $\beta$ is a function which assigns to each
border node $n$ a subset $\beta (n) $ of $V(G)$; and 
\item $\gamma$ is an injection from $N(\Gamma)$ into $V(G)$
\end{itemize}
satisfying the axioms below.  As before, for
each $X \subseteq N(\Gamma)$ we denote
$\{  \gamma (n) : n \in X\}$ by
$\gamma (X)$, and for each $c_0 \in C (\Gamma )$ we denote
$\cup ( \alpha (c) : c \in C(\Gamma) - \{c_0\})$ by
$\alpha (- c_0)$.
For each cuff $\Theta$ let $\alpha(\Theta)$ be the union of $\alpha(c)$ over all cells that border $\Theta$.
The axioms are as follows:

\begin{itemize}
\item[{\bf (J1)}] $G = \bigcup (\alpha (c) : c \in C( \Gamma))$, and $E(\alpha (c) \cap \alpha (c')) = \emptyset$ for distinct cells $c, c'$;

\item[{\bf (J2)}] $\beta(n)\cap N(\Gamma) = \emptyset$ for each $n \in \bd (\Sigma )$, and $\beta(n) \subseteq V ( \alpha (c))$ 
for each border cell $c$ and each $n \in \tilde{c}$;

\item[{\bf (J3)}] If $c$ is an internal cell then
$V(\alpha (c) \cap \alpha (-c)) = \gamma ( \tilde{c})$;

\item[{\bf (J4)}] If $c$ is a border cell and $\tilde{c}=\{n_1, n_2\}$, then
$$ V( \alpha (c) \cap \alpha (-c)) = \gamma ( \tilde{c})\cup \beta (n_1) \cup \beta(n_2); $$

\item[{\bf (J5)}] If $n_1,n_2,n_3,n_4$ are nodes bordering the same cuff of $\Sigma$ and in order, then
$\beta (n_1) \cap \beta (n_3) \subseteq \beta(n_2) \cup \beta (n_4)$;

\item[{\bf (J6)}] For each cuff $\Theta$, all the sets $\beta(n)\;(n\in N(\Gamma)\cap \Theta)$ have the same cardinality, called the {\em depth}
of $\Theta$ in the portrayal.

\item[{\bf (J7)}] For each cuff $\Theta$, there is a 3-ring $F_{\Theta}$ around $\Theta$, and for all distinct cuffs, the corresponding annuli 
are disjoint. 

\item[{\bf (J8)}] For each internal cell $c$, and all $n_1,n_2\in \tilde{c}$, there is a path of $\alpha(c)$ between $\gamma(n_1), \gamma(n_2)$
using no other vertex in $\gamma(\tilde{c})$. For each border cell $c$, with $\tilde{c}=\{n_1,n_2\}$, either there is a 
path of $\alpha(c)$ between $\gamma(n_1), \gamma(n_2)$
using no vertices of $\beta(n_1)\cup \beta(n_2)$, or there is an internal cell incident with both $n_1,n_2$. 

\end{itemize}
The {\em depth} of the portrayal is the maximum of the depths of its cuffs (or zero, if there are no cuffs).

Let $(\Sigma, \Gamma, \alpha,\beta,\gamma)$ be a clean portrayal of $G$, and let $\mac T$ be a tangle in $G$ of order $\theta$. 
We say:
\begin{itemize}
\item
$(\Sigma, \Gamma, \alpha,\beta,\gamma)$ is {\em $\mac T$-central} if
$(\alpha(c),\alpha(-c))\in \mac T$ for each cell $c$ of $\Gamma$;
\item $\mac T$ is {\em matted} if $A$ is planar for each $(A,B)\in \mac T$;
\item 
$\mac T$ is {\em $z$-respectful} if for every $\Gamma$-normal O-arc in $\hat{\Sigma}$ with $|F\cap N(\Gamma)|\le z$, there is a closed disc $D$ in $\hat{\Sigma}$
with boundary $F$, such that $(\alpha(D), \alpha(\hat{\Sigma}\setminus D)\in \mac T$, and if $F\subseteq \Sigma$ then $D\subseteq \Sigma$. We write $\ins(F)=D$.
\end{itemize}

Let $(\Sigma, \Gamma, \alpha,\beta,\gamma)$ be a clean portrayal of $G$, and let $\sigma$ be a side of a cell, with ends 
$n_1,n_2$. If there is a path of $\alpha(c)$ between $\gamma(n_1), \gamma(n_2)$
using no other vertex in $\gamma(\tilde{c})$, and with no other vertex in $\beta(n_1)\cup \beta(n_2)$ if $c$ is a border cell,
let us fix one such path, which we call the {\em support} for $\sigma$. If $c$ has two sides
$\sigma,\sigma'$ we choose the corresponding supports to be equal. If $c$ has three sides, 
let us choose the corresponding supports such that their union is either a cycle or a tree (which therefore has exactly one vertex of degree three).
In this section we prove:

\begin{thm}\label{cleanp}
For every connected surface $\Sigma$ with $q$ cuffs, there is a function $f$ with the following property.
Let $d,k\ge 0$ be integers, and
let $\mac T$ be a matted tangle in a graph $G$.
Suppose that $G$ admits a $\mac T$-central clean portrayal of depth $d$ in $\Sigma$, such that $\mac T$ is 6-respectful; and 
suppose that $\mac T$ has order $\theta\ge f(k,d)$.
Then either:
\begin{itemize}
\item $G$ has K-number at least $k$; or
\item there exists $Z\subseteq V(G)$ with $|Z|\le f(k,d)$ such that $G\setminus Z$ can be drawn in $\Sigma$.
\end{itemize}
\end{thm}
\Proof
Let $g$ be the function of \ref{3planar2} for the surface $\hat{\Sigma}$, and define $f(k,d) = 2k(d+3)q+g(k)+4$. 
We will show that $f$ satisfies the theorem.
Let $(\Sigma, \Gamma, \alpha,\beta,\gamma)$ be a clean portrayal of $G$. For each cuff $\Theta$, let $d_{\Theta}$ be its depth.
Our first objective is to prove that for every cuff $\Theta$, 
we can delete a few vertices from a small neighbourhood of the cuff, and rearrange the remainder  of the portrayal,
to become a portrayal on the surface where the cuff $\Theta$ has been capped. Thus, suppose that $\Theta$ is a cuff of $\Sigma$. Let
$\Sigma_{\Theta}$ be the annulus bounded by $\Theta$  and $F_{\Theta}$. 
Let the nodes of $F_{\Theta}$ be $n_1\CC n_t\DD n_1 $ in order.
Let $(H_{\Theta},\Omega_{\Theta})$ be the society where 
$H_{\Theta}=\alpha(\Sigma_{\Theta})$, and $\Omega_{\Theta}$ is the permutation
$\gamma(n_1)\rightarrow \gamma(n_2)\rightarrow\cdots\rightarrow  \gamma(n_t)\rightarrow \gamma(n_1)$.

For each node $n\in F_{\Theta}$, there is a $\Gamma$-normal line segment $L(n)$ with one end $n$ and the other some node $p(n)\in \Theta$, 
containing exactly three nodes.
Let us choose $L(n)\;(n\in N(\Gamma)\cap F_{\Theta})$ such that for all nodes $n_1,n_2\in F_{\Theta}$, either $L(n_1), L(n_2)$ are disjoint,
or $p(n_1)= p(n_2)$ and $L(n_1)\cap L(n_2)=\{p(n_1)\}$, or $p(n_1)= p(n_2)$ and $L(n_1)\cap L(n_2)$ is a line segment between $p(n_1)$ 
and some other node of $\Gamma$.
For $1\le i\le t$, let $\sigma_i$ be the interior of the line segment with
ends $n_i,n_{i+1}$ included in $F_{\Theta}$ (reading subscripts modulo $t$).
For $1\le i\le t$, let $\phi_i$ be the closed curve obtained by following $L(n_i)$ from $p(n_i)$ to $n_i$, following $\sigma_i$ to $n_{i+1}$,
following $L(n_{i+1})$ to $p(n_{i+1})$, and then following part of $\Theta$ back to $p(n_i)$, in such a way that $\phi_i$
is null-homotopic in $\Sigma_{\Theta}$ in the natural sense. (Since $L(n_i), L(n_{i+1})$ may have two nodes in common, $\phi_i$
does not necessarily trace an O-arc.) If we cut $\Sigma_{\Theta}$ along $L(n_i)$ and then along $L(n_{i+1})$, we obtain exactly two 
closed discs. We define $D_{i,i+1}$ to be the one of these discs that includes $\sigma_i$. For $0\le i<j\le t$ let $D(i,j)$
be the union of $D_{i,i+1},D_{i+1,i+2}\LL D_{j-1,j}$. 
\\
\\
(1) {\em Not all the nodes $p(n)\;(n\in N(\Gamma)\cap F_{\Theta})$ are equal.}
\\
\\
Suppose all the $p(n)$ are equal to $m$ say. 
Since each $\phi_i$ is null-homotopic in $\Sigma_{\Theta}$, it follows that 
their concatenation (in order) 
is also null-homotopic
in $\Sigma_{\Theta}$; that is, the curve obtained by following $L(n_t)$ from $m$ to $n_t$, following $F_{\Theta}$ along its entire length 
and thereby reaching 
$n_t$ again, and then following $L(n_t)$ from $n_t$ to $m$, is null-homotopic in the annulus, a contradiction. This proves (1).

\bigskip

From (1), we may assume that $p(n_0)\ne p(n_s)$ for some $s\in \{1\LL t\}$, and consequently $L(n_0)$ and $L(n_s)$ are disjoint.
If we cut $\Sigma_{\Theta}$ along $L(n_0)$ and $L(n_s)$, we obtain two discs $D_1,D_2$, where $n_0\LL n_s\in D_1$
and $n_s,n_{s+1}\LL n_t\in D_2$.
If $0\le i\le s$ then $L_i\subseteq D_1$, and otherwise $L_i\subseteq D_2$. 
For $0\le q\le r\le s$, some of $L(n_1), L(n_q), L(n_r), L(n_s)$ may intersect one another; nevertheless, if we cut $D_1$
along the portions of $L(n_q)$ and $L(n_r)$ that are in the interior of $D_1$, we obtain exactly three discs. Exactly one
of them is bounded by a portion of $L_q$ and a portion of $L_r$, and we denote it by $D(q,r)$. Let $H(q,r)=\alpha(D(q,r))$,
and let $\Omega(q,r)$ be the cyclic permutation 
$$\gamma(n_q)\rightarrow\gamma(n_{q+1})\rightarrow\cdots\rightarrow\gamma(n_r)\rightarrow\gamma(n_q).$$

For $0\le q<s$, its {\em successor}, if it exists, is the smallest $r$ with $q<r\le s$ such that there is a cross of 
$(H(q,r),\Omega(q,r))$ 
vertex-disjoint from $\gamma(N(\Gamma)\cap L(n_q))\cup \beta(p(n_q))$. 
\\
\\
(2) {\em Suppose that $0\le q\le r\le s$ and $r$ is the successor of $q$. Let $W$ be the union of $\alpha(c)$ over all 
cells $c$ of $\Gamma$ with $c\subseteq D(q,r)$, together with the supports of all sides $\sigma\subseteq \bd (D(q,r))$ 
of cells of $\Gamma$ such that $\sigma$ has a support.
Then there is a K-graph in $W$ containing no vertices in $\gamma(N(\Gamma)\cap L(n_q))\cup \beta(p(n_q))$.}
\\
\\
By \ref{getnest}, there is a 3-nest $(R_1,R_2,R_3)$
where $R_1=\Theta$ and $R_3=F_{\Theta}$.
Since $r$ is the successor of $q$, it follows that
$(H(q,r),\Omega(q,r))$ has a cross with components $Q,R$ disjoint from $\gamma(N(\Gamma)\cap L(n_q))\cup \beta(p(n_q))$, and 
so $p(n_q)\ne p(n_r)$. Thus the boundary of $D(q,r)$ is the union of four line segments $L(n_q)$, $D(q,r)\cap R_3$,
$L(n,r)$ and $D(q,r)\cap R_1$. Let $t=r-q+5$, and let the nodes in $\bd(D(p,q))$ be 
be $m_1\LL m_S\LL m_T$ in order, where $m_1=p(n_q)$, $m_i = n_{q+i-3}$ for $3\le i\le r-q+3$, and $m_S=p(n_r)$, and the nodes in
$R_1\cap D(p,q)$ are $m_1,m_T,m_{T-1}\LL m_S$ in order.

Let $A$ be the union of $\alpha(c)$ over all border cells $c$ of $\Gamma$ with $c\subseteq D(q,r)$. 
Let $B$ be the union of $\alpha(c)$ over all internal cells $c$ with $c\subseteq D(p,q)$, together with the supports of all sides $\sigma\subseteq \bd (D(q,r))$ 
of cells of $\Gamma$ that are not subsets of $R_1$. Thus $(A,B)$ is a separation of $W$,
with $A\cap B=\gamma(N(\Gamma)\cap R_1)$. 
From the way 
we chose supports, it follows that $B$ can be drawn in the disc $D(q,r)$ in such a way
that $\gamma(n)$ is drawn at $n$ for each node $n\in \bd (D(q,r))$.

Now $Q,R$ will contain vertices in $A\cap B$, but we can choose $Q,R$
such that every subpath $S$ of one of them, maximal with no internal vertex in $\gamma(N(\Gamma))$, is either a path of $A$
or a support included in $B$;
and consequently $Q,R$ are paths of $A\cup B$.

For $i = 2,3$, let $P_i$ be the union of the supports of the sides of $\Gamma$ included in $R_i\cap D(q,r)$.  But $P_1$ is more complicated,
because the sides of cells in $R_1$ need not have supports. If $\sigma\subseteq R_1$ is a side with a support, take such a support;
and otherwise there is a side of an internal cell incident with the same two nodes, and take the support of that side. Let $P_1$
be the union of all the supports chosen. Thus, for $i = 1,2,3$, $P_i$ is a path
of $W$ between $\gamma(m_i)$ and $\gamma(m_{t+1-i})$ for $i = 1,2,3$. Also, $P_1$ is a path of $A\cup B$, and $P_2,P_3$ are paths of $B$, 
and these three paths are pairwise vertex-disjoint. 
From \ref{crosstoK}, the union of the five paths $P_1,P_2,P_3,Q,R$ is not planar.
Since neither of $Q,R$ contain any of $\gamma(m_1), \gamma(m_2), \gamma(m_3)$, it follows that these three vertices have degree one in the union of the five paths; and so there is a K-graph $G_j$ in the union of the five paths that contains none of 
of $\gamma(m_1), \gamma(m_2), \gamma(m_3)$. Moreover, it contains no vertices in $\beta(m_1)$ since such vertices belong 
to none of the five paths. This proves (2).

\bigskip

Now there are two cases: perhaps $0$ has a successor
$i_1$, and $i_1$ has a successor $i_2$, and so on at least $k$ times; or this does not happen. 
\\
\\
(3) {\em If there exist $0=i_0<i_1<i_2<\cdots<i_k\le s$ such that $i_j$ is the successor of $i_{j-1}$ for $1\le j\le k$,
then the K-number of $G$ is at least $k$.}
\\
\\
From (2), for $0\le j\le k-1$ there is a K-graph as in (2); and we claim these K-graphs are pairwise vertex-disjoint. Indeed,
let $0\le q<r\le q'<r'\le s$ and suppose that $r$ is the successor of $q$ and $r'$ is the successor of $q'$. Let $W$ be the 
union of $\alpha(D(q,r))$
with the supports of all sides $\sigma\subseteq \bd (D(q,r))$
of cells of $\Gamma$; and define $W'$ similarly for $q', r'$. We claim that $V(W\cap W')=\emptyset$. For suppose that
$v$ belongs to both subgraphs. If $v=\gamma(n)$ for some node $n$, then $n$ is incident with a cell $c\subseteq D(q,r)$
and incident with a cell $c'\subseteq D(q',r')$. Consequently $n$ belongs to the intersection of the boundaries of 
$D(q,r)$ and $D(q',r')$, and hence to $L(n_{q'})$, contrary to the choice of $W'$. So we assume that 
$v\notin \gamma(N(\Gamma))$. If $v\in V(\alpha(c))$ for some internal cell $c$, then $c$ is unique, and so there is a side of $c$
included in $D(q,r)$ and a side of $c$ included in $D(q',r')$, a contradiction to the definition of a $k$-ring.
Thus $v\in V(\alpha(c))$ for some border cell $c$ included in $D(q,r)$, and $v\in V(\alpha(c'))$ for some border cell $c'$ 
included in $D(q',r')$. But then $v\in \beta(n_{q'})\cup \beta(n_q)$, again a contradiction. This proves (3).

\bigskip

We may therefore assume that there do not exist $i_0\LL i_k$ as in (3). Choose $h$ maximum such that there exist
$0=i_0<i_1<i_2<\cdots<i_h\le s$ such that $i_j$ is the successor of $i_{j-1}$ for $1\le j\le h$. Thus $i_h$ has no successor.
Let $X_1$ be the union of the $h$ sets $\gamma(L(n_{i_j})\cup \beta(n_{i_j})\;(0\le j\le h)$, and choose $X_2$ similarly 
for $D_2$. Let 
$$X_{\Theta}=X_1\cup X_2\cup \gamma(L(n_{1})\cup \beta(n_{1}) \cup \gamma(L(n_{s})\cup \beta(n_{s});$$
then every cross of $(H_{\Theta},\Omega_{\Theta})$ has a vertex in $X_{\Theta}$, and so $(H_{\Theta},\Omega_{\Theta})$ has urbanity at most $2k(d_{\Theta}+3)$.
Thus we have proved:
\\
\\
(4) {\em There exists $X_{\Theta}\subseteq \gamma(N(\Gamma)\cap \Sigma_{\Theta})$ with $|X_{\Theta}|\le 2k(d_{\Theta}+3)$
meeting every cross of $(H_{\Theta},\Omega_{\Theta})$.}

\bigskip
Next, we need:
\\
\\
(5) {\em There is no separation $(A,B)\in \mathcal{T}$ with $G\subseteq A\cup \alpha(\Sigma_{\Theta})$ of order 
less than $\theta-6-2d_{\Theta}$.}
\\
\\
Suppose that such $(A,B)$ exists. Let $D(i,i+1)\;(0\le i\le t)$ be as before, let 
$D_i= D(0,1)\cup D(1,2)\cupcup D(i-1,i)$ for $0\le i\le t$ (so $D(0) = \emptyset$), let $G_i=\alpha(D_i)$,
and let $B_i=B\cap G_i$. Let $H_i$ be the union of $\alpha(c)$
over all cells $c$ with $c\not\subseteq D_i$ (so $G_i\cup H_i=G$), and let $A_i=A\cup H_i$. Thus $(A_i,B_i)$ is a separation of $G$. 
We claim that 
$$V(A_i\cap B_i)\subseteq V(A\cap B)\cup \gamma(L_1)\cup \beta(p(n_1))\cup \gamma(L_i)\cup \beta(p(n_i)).$$
To see this, suppose that $v\in V(A_i\cap B_i)$ belongs to none of the five sets on the right side above.
It 
follows that $v\in V(B)\setminus V(A)$, and $v\in V(G_i\cap H_i)$. The latter implies that either $v=\gamma(n)$ for some node 
$n$ in the boundary of $D_i$, since $v\notin \beta(p(n_1))$ and $v\notin \beta(p(n_i))$. Since $v\notin \gamma(L_1)\cup  \gamma(L_i)$, 
it follows that either $n\in \Theta\setminus \{p(n_1),p(n_i)\}$, or
$n\in F_{\Theta}\setminus \{n_1,n_i\}$.
It is impossible that $n\in \Theta\setminus \{p(n_1),p(n_i)\}$, since $v\in H_i$. If $n\in F_{\Theta}\setminus \{n_1,n_i\}$, 
then since $v\in A\cup H_i$, it follows that $v\in V(A)$ 
(because $\alpha(c)\subseteq A$ for every cell $c$ incident with $n$ that is not included in $D_i$), a contradiction. 
This proves our claim, and we deduce that the order of $(A_i,B_i)$ is at most $|V(A\cap B)|+6+2d_{\Theta}$, and so one of $(A_i,B_i),(B_i,A_i)\in \mac T$.
Since $(B_0,A_0)\in \mac T$ (because $B_0$ is the null graph),
we may choose $i$ with $0\le i\le t$ maximum such that $(B_i,A_i)\in \mac T$; and since
$(B_t,A_t)\notin \mac T$ (since $B_t=B\cap \alpha(\Sigma_{\Theta})$ and so $A\cup B_t=G$)
it follows that $i<t$.
Consequently $(A_{i+1},B_{i+1})\in \mac T$, since it has order at most 
$|V(A\cap B)|+6+2d_{\Theta}< \theta$. 
Let $P=\alpha(D(i,i+1)$ and let $Q$ be the union of $\alpha(c)$ over all cells $c$ not 
included in $D(i,i+1)$. Then $(P,Q)$ is a separation of order at most $6+2d_{\Theta}$, and belongs to $\mac T$ since $\mac T$
is 6-respectful.
But $P\cup B_i\cup A_{i+1}=G$, contradicting the third tangle axiom. This proves (5).

\bigskip

Let us paste a disc onto $\Theta$, forming $\Sigma'$ say with one fewer cuff. Let $\Delta_{\Sigma}$ be the disc in 
$\Sigma'$ bounded by $F_{\Theta}$ that includes the disc we just added.
Since $H_{\Theta}\setminus X_{\Theta}, \Omega_{\Theta}\setminus X_{\Theta})$ has no cross, we can apply \ref{nocross} to it, 
and deduce that it admits a rural portrayal $(\Sigma_1, \Gamma_1, \alpha_1, \gamma_1)$. in a sphere $\Sigma_1$; 
and we may assume that $\Delta_{\Theta}$ is the disc of $\Sigma_1$ including $U(\Gamma_1)$, and for each $n\in F_{\Theta}$, 
if $\gamma(n)\notin X_{\Theta}$
then $\gamma_1(n)=\gamma(n)$. 
Consequently, 
we can obtain a portrayal $(\Sigma', \Gamma', \alpha', \beta', \gamma')$ of $G\setminus X_{\Theta}$ in $\Sigma'$ 
by combining parts of the two portrayals in the natural way. 
\\
\\
(6) {\em 
There exists $X\subseteq V(G)$ with $|X|\le 2k(d+3)q$, and a matte portrait 
$(\Gamma_0, \alpha_0, \gamma_0)$ of $G\setminus X$ in $\hat{\Sigma}$, such that for each cell $c_0$ of $\Gamma_0$, either
there is a cell $c$ of $\Gamma$ with $\alpha_0(c_0)\subseteq \alpha(c)$, or $\alpha(c_0)\subseteq \alpha(\Sigma_{\Theta})$ for some cuff $\Theta$ of $\Sigma$.}
\\
\\
For each cuff $\Theta$ of $\Sigma$, $|X_{\Theta}|\le d_{\Theta}+3$. 
By repeating the process just described for each cuff of $\Sigma$, we obtain $X\subseteq V(G)$ with $|X|\le 2k(d+3)q$, and a portrayal
$(\hat{\Sigma}, \Gamma_0, \alpha_0, \gamma_0)$ of $G\setminus X$. We need to show that it is matte. Let $c_0$ be a cell of $\Gamma_0$.
Choose a separation $(A,B)$ of $G$ with $X\subseteq V(A\cap B)$, such that $B\setminus X=\alpha_0(c_0)$. Since $\mac T$
has order at least $|X|+3$, one of $(A,B),(B,A)\in \mac T$. 
If $c_0$ is not a subset of $\Delta_{\Theta}$ for any cuff $\Theta$, then there is a cell $c$ of $\Gamma$ with $c_0\subseteq c$
and $\alpha_0(c_0)\subseteq \alpha(c)$; and consequently $\alpha_0(c_0)$ is planar, since $\alpha(c)$ is planar (because
$(\alpha(c),\alpha(-c))\in \mathcal{T}$ and $\mathcal{T}$ is matted). Thus we may assume that $c_0\subseteq \Delta_{\Theta}$
for some cuff $\Theta$ of $\Sigma$. Hence $A\cup \alpha(\Sigma_{\Theta})=G$, and so $(A,B)\notin \mac T$ by (5). Hence $(B,A)\in \mac T$, and so $B$ is planar, since $\mac T$ is matted; and hence $\alpha_0(c_0)$ is planar, and therefore 
$(\hat{\Sigma}, \Gamma_0, \alpha_0, \gamma_0)$ is matted. This proves (6).

\bigskip

From \ref{3planar2} applied to $G\setminus X$ and the portrait $(\Gamma_0, \alpha_0, \gamma_0)$ in $\hat{\Sigma}$, we deduce that either:
\begin{itemize}
\item $G\setminus X$ has K-number at least $k$; or
\item there exists $Z\subseteq V(G\setminus X)$ with $|Z|\le g(k)$ such that $(G\setminus X)\setminus Z$ can be drawn in $\hat{\Sigma}$; or
\item there is a separation $(A',B')$ of $G\setminus X$ of order at most $g(k)$ such that both $A',B'$ are not planar.
\end{itemize}
In the first and second cases the theorem holds, since $|X|+|Z|\le f(k,d)$. Suppose $(A',B')$ is as in the third case.

Choose a separation $(A,B)$ of $G$ with $X\subseteq V(A\cap B)$ and $A\setminus X=A'$ and $B\setminus X=B'$. Since one of 
$(A,B),(B,A)\in \mac T$ and $\mac T$ is matted, it follows that one of $A,B$ is planar, a contradiction. So the third case does not hold.
This proves \ref{cleanp}.~\bbox

\section{The proof of \ref{thm1}}

We need the following:
\begin{thm}\label{oldtonew}
For every graph $H$, there are integers $p,q,z\ge 0$ and $\theta>z$, such that, for every graph $G$ with no $H$ minor,
and every matted tangle in $G$ of order at least $\theta$, there exists $Z\subseteq V(G)$ with $|Z|\le z$ 
and a clean $(\mathcal{T} \setminus Z)$-central portrayal
$\Pi = (\Sigma, \Gamma, \alpha, \beta, \gamma)$ of  $G \setminus Z$ with depth
$\le p$, such that $\Sigma $ has $\le q$ cuffs,
and $H$ cannot be drawn in $\Sigma$.
\end{thm}
\Proof
Let $H$ be the graph consisting of $k$ disjoint copies of $K_5$. For each surface $\Sigma$
in which $H$ cannot be drawn, and all integers $p,p' \ge 0$, define $\sigma(\Sigma,p,p') = 4p+5$.
The function $\sigma$ is what is called a ``standard'' in~\cite{GM17}.
From theorem (13.4) of~\cite{GM17}, there are integers $p,q,z\ge 0$ and $\theta>z$, such that, for every graph with no $H$ minor, 
and every tangle in $G$ of order at least $\theta$, there exists $Z\subseteq V(G)$ with $|Z|\le z$
and a $(\mathcal{T} \setminus Z)$-central portrayal
$\Pi = (\Sigma, \Gamma, \alpha, \beta, \gamma)$ of  $G \setminus Z$ with depth
$\le p$, such that $\Sigma $ has at most $q$ cuffs,
$H$ cannot be drawn in $\Sigma$, and $\Pi$ is
$(4p+5)$-redundant and true.
We will avoid defining ``portrayal'' or ``redundant'' or ``true'' here, and instead rephrase this result in the language of this paper.
We may assume that $z\ge 6$ and $\theta> z+4p+6$, and so $\mac T\setminus Z$ has order more than $4p+6$.
A portrayal in the sense of~\cite{GM17} is again a 5-tuple $(\Sigma, \Gamma, \alpha,\beta,\gamma)$, and it will turn out to be a clean portrayal. We see
this as follows. The conditions on the painting $\Gamma$ are same for portrayals and clean portrayals, because of theorems (8.2)
and (8.4) of~\cite{GM17}; and portrayals satisfy
{\bf (J1)} and {\bf (J5)}. Because the portrayal is $(4p+5)$-redundant and
$\theta> 4p+5$, theorems (9.1) and (9.8) of~\cite{GM17} imply that {\bf (J8)} holds.
Theorem (9.8) of~\cite{GM17} implies {\bf (J6)}, and theorem (9.7) of~\cite{GM17} implies {\bf (J2)}. Theorem (8.7) implies {\bf (J3)}
and {\bf (J4)}. Since $\theta> z+4p+6$,
theorems (6.3) and (6.4) of~\cite{GM17} imply that $\mac T\setminus Z$ is $6$-respectful.

We still need to verify {\bf (J7)}.
Every region of $\Gamma$ (regarded as a painting in $\hat{\Sigma}$) is an open disc, by theorem (8.1) of~\cite{GM17},
and so there is a ``radial drawing'' $\Gamma'$ in the sense of~\cite{GM11}.
Since $\mac T\setminus Z$ is $z$-respecful,  there is a ``slope'' in $\Gamma'$
of order $z$. As in~\cite{GM11}, we may speak of ``restraints''  in $\Gamma'$. By theorem (8.10) of~\cite{GM11},
for each cuff $\Theta$
of $\Sigma$, which therefore bounds a region $r$ say of $\hat{\Sigma}$, the union of all restraints that include $r$  with ``length''
at most $2z$ is simply-connected, and does not meet any other cuff of $\Sigma$ (because of theorems (6.3) and (6.4) of~\cite{GM17}).
Consequently we can satisfy {\bf (J7)} within the union of the insides of these restraints.
This proves \ref{oldtonew}.~\bbox

Now we can complete the proof of \ref{thm1}, which we restate:
\begin{thm}\label{thm1again}
For all $k\ge 0$, there is a number $f(k)$, such that for every graph $G$
with K-number at most $k$, there exists $X\subseteq V(G)$ with $|X|\le f(k)$ such that $G\setminus X$ can be drawn in 
a surface with genus at most $k$.
\end{thm}
\Proof  We proceed by induction on $k$. The result is true by the Kuratowski-Wagner theorem if $k=0$, so we assume that $k>0$, and we may assume that $f(i)\le f(k-1)$ for $0\le i\le k-1$.
Let $H$ be the graph consisting of $k+1$ disjoint copies of $K_5$. Choose $p,q\ge 0$ and $z\ge 6$ and $\theta>z+4p+6$ 
satisfying \ref{oldtonew}.
$f(k)/3\le f(k)-2f(k-1)$
For each connected surface $\Sigma$ with at most $q$ cuffs, and with genus at most $k$, let $f_{\Sigma}$ be the function $f$ as in \ref{cleanp}.
Choose a number $F$ such that $f_{\Sigma}(k,d)\le F$ for all such $\Sigma$. By increasing $\theta$, we may assume that 
$\theta> z+2(k+1)(d+3)q$.
Let $f(k)$ satisfy 
$$f(k)>\max(3f(k-1), 2f(k-1)+ \theta, z+F).$$
We will prove that 
$f(k)$ satisfies the theorem.

Let $G$ be a graph with K-number at most $k$; we must show that there exists $X\subseteq V(G)$ with $|X|\le f(k)$ such that 
$G\setminus X$ can be drawn in a surface of genus at most $k$.  Let $\mac T$
be the set of all separations $(A,B)$ of $G$ of order less than $\theta$ such that $A$ is planar.
\\
\\
(1) {\em We may assume that $\mathcal{T}$ is a matted tangle in $G$
of order $\theta$.}
\\
\\
Suppose first that there is a separation $(A,B)$ of $G$ with order
less than $\theta$
such that both $A,B$ are nonplanar. Let $Y=V(A\cap B)$, and let the K-numbers of $A\setminus Y$, $B\setminus Y$ 
be $k_1,k_2$ respectively. Since $A$ is nonplanar, it follows that $k_2<k$, and similarly $k_1<k$; and $k_1+k_2\le k$.
From the inductive hypothesis, there exist $X_1\subseteq V(A\setminus Y)$ with $|X_1|\le f(k_1)$ such that 
$A\setminus V(A\cap B)\setminus X_1$ can be drawn in a surface $\Sigma_1$ of genus at most $k_1$. Define $X_2, \Sigma_2$ similarly for 
$B\setminus Y$. Let
$X=X_1\cup X_2\cup Y$; then $|X|\le f(k_1)+f(k_2)+ \theta\le f(k)$, and $G\setminus X$ can be drawn in the disjoint union of $\Sigma_1,\Sigma_2$, which has genus at most $k_1+k_2\le k$.

So we may assume that there is no such $(A,B)$. 
It follows that 
for every separation $(A,B)$ of $G$ of order less than $\theta$, one of $(A,B),(B,A)\in \mac T$.
Suppose that $(A_i,B_i)\in \mac T$ for $i = 1,2,3$, and let $X_i=V(A_i\cap B_i)$ for $ 1\le i\le 3$. 
If $G\setminus (X_1\cup X_2\cup X_3)$ is planar, then the theorem holds, since $f(k)\ge 3(\theta-1)$. So we may assume that
there is a K-graph $H$ in $G$ disjoint from $X_1\cup X_2\cup X_3$. For $i = 1,2,3$, since $A_i$ is planar it follows that
$H\not\subseteq A_i$; and since $H$ is connected and $(A_i,B_i)$ is a separation and $V(H)\cap X_i=\emptyset$, it follows
that $H, A_i$ are vertex-disjoint. Hence $A_1\cup A_2\cup A_3\ne G$. It follows that $\mathcal{T}$ is a matted tangle in $G$
of order $\theta$. This proves (1).

\bigskip

From \ref{oldtonew},
there exists $Z\subseteq V(G)$ with $|Z|\le z$
and a clean $(\mathcal{T} \setminus Z)$-central portrayal
$(\Sigma, \Gamma, \alpha, \beta, \gamma)$ of  $G \setminus Z$ with depth
$\le p$, such that $\Sigma $ has $\le q$ cuffs,
and $H$ cannot be drawn in $\Sigma$. 
Since $H$ cannot be drawn in $\hat{\Sigma}$, 
it follows that $\hat{\Sigma}$ has genus at most $k$.
By \ref{cleanp}, applied to $G\setminus Z$ (with $k$ replaced by $k+1$), and since $\mac T\setminus Z$ has order at least
$\theta-z\ge 2(k+1)(d+3)q$,
and $G\setminus Z$ has $K$-number at most $k$, it follows that 
there exists $Z'\subseteq V(G)\setminus Z$ with $|Z'|\le f_{\Sigma}(k+1,p)\le F$ such that $G\setminus (Z\cup Z')$ can be drawn in $\Sigma$.
Since $|Z\cup Z'|\le z+F\le f(k)$, this proves \ref{thm1}.~\bbox



\begin{thebibliography}{99}
\bibitem{archdeacon} D. Archdeacon, ``A Kuratowski theorem for the projective plane'', 
{\em J. Graph Theory}, {\bf 5} (1981), 243--246.

\bibitem{battle} J, Battle,  F. Harary, Y. Kodama andi J. W. T. Youngs, ``Additively of the genus of a
graph'', {\em Bull. Amer. Math. Soc.} {\bf 68} (1962), 565--568. 

\bibitem{bienstock} D. Bienstock and N. Dean, ``On obstructions to small face covers in planar graphs'',
{\em J. Combinatorial Theory, Ser. B}, {\bf 55} (1992), 163--189.

\bibitem{ding} G. Ding, P. Seymour and P. Winkler, ``Bounding the vertex cover number of a
hypergraph'', {\em Combinatorica,} {\bf 14} (1994), 23--34.


\bibitem{kuratowski}  K. Kuratowski,
``Sur le probl\`{e}me des courbes gauches en topologie'',
{\em Fund. Math.} 15 (1930), 271--283.

\bibitem{miller} G. L. Miller, ``An additivity theorem for the genus of a graph '', {\em J. Combinatorial Theory, Ser. B}, 
{\bf 43} (1987), 25--47.

\bibitem{myrvold} W. Myrvold and J. Woodcock, ``A large set of toroidal obstructions and how they were discovered'',
{\em Electronic J. Combinatorics} {\bf 25} (2018), \#P1.16.

\bibitem{Ksums} N. Robertson and P. Seymour, ``Excluding sums of Kuratowski graphs'', in preparation.

\bibitem{GM5} N. Robertson and P. Seymour,  ``Graph minors. V. Excluding a planar graph'',
{\em J. Combinatorial Theory, Ser. B,} {\bf 41} (1986), 92--114.

\bibitem{GM8} N. Robertson and P. Seymour, ``Graph minors. VIII. A Kuratowski theorem for
general surfaces'', {\em J. Combinatorial Theory, Ser. B,} {\bf 48} (1990), 255--288.

\bibitem{GM9} N. Robertson and P. Seymour, ``Graph minors. IX. Disjoint crossed paths'', {\em J.
Combinatorial Theory, Ser. B,} {\bf 49} (1990), 40--77.

\bibitem{GM10} N. Robertson and P. Seymour, ``Graph minors. X. Obstructions to tree-decomposition'',
{\em J. Combinatorial Theory, Ser. B,} {\bf 52} (1991), 153--190.

\bibitem{GM11} N. Robertson and P. Seymour, ``Graph minors. XI. Circuits on a surface'', {\em J. Combinatorial
Theory, Ser. B}, {\bf 60} (1994), 72--106.

\bibitem{GM16} N. Robertson and P. Seymour, ``Graph minors. XVI. Excluding a non-planar graph'',
{\em J. Combinatorial Theory, Ser. B}, {\bf 89} (2003), 43--76.

\bibitem{GM17} N. Robertson and P. Seymour, ``Graph minors. XVII.  Taming a vortex'',
{\em J. Combinatorial Theory, Ser. B}, {\bf 77} (1999), 162--210.

\bibitem{GM20} N. Robertson and P. Seymour, ``Graph minors. XX. Wagner's conjecture'', {\em J.
Combinatorial Theory, Ser. B}, {\bf 92} (2004), 325--357.

\bibitem{wagner} K. Wagner,
``\"Uber eine Eigenschaft der ebenen Komplexe'',
{\em Math. Ann.} 114 (1937), 570--590.

\end{thebibliography}
\end{document}